\documentclass[a4paper,10pts]{article}
\makeatletter
   \let\accent@spacefactor\relax
 \makeatother       
\usepackage[french]{babel}
\usepackage{amsmath,amscd,amssymb}
\usepackage{fancyheadings}

\def \p {{\bf P}}

\def \pu {{\bf P}^1}

\def \n2 {{n+1\over{2}}}
\def \g {{\bf G}}

\def \fl {\longrightarrow}

\def \vs {\vskip}

\def \I {{\cal I}}

\def \oo {{\cal O}}

\def \dm {{\textit{D\'emonstration}}}
\def \gh {{\mathfrak h}}
\def \gg {{\mathfrak g}}
\def \gb {{\mathfrak b}}
\def \gp {{\mathfrak p}}
\def \gq {{\mathfrak q}}

\def \H {{\bf Hilb}}

\def \a {{\alpha}}
\def \b {{\beta}}
\def \gsp {{G/P}}
\def \vp {\varphi}
\def \opu {{{\cal O}_{\pu}}}
\def \S {{\Sigma}}
\def \pic {{\rm{Pic}}}

\def \noi {\noindent}
\def \Mor #1#2{{\bf{Hom}}_{#1}(\pu,#2)}
\def \mapdown#1{\Big\downarrow
\rlap{$\vcenter{\hbox{$\scriptstyle#1$}}$}}
\def \G {{\bf G}}

\begin{document}

{\centerline {\Large {\bf{Courbes rationnelles sur les
vari\'et\'es homog\`enes}}}}
{\centerline {\Large {\bf{et une d\'esingularisation plus
fine des vari\'et\'es de Schubert}}}}

\vskip 0.7 cm

{\centerline { Nicolas PERRIN }}
{\centerline {Universit\'e de Versailles St-Quentin }}
{\centerline {45 Avenue des Etats-Unis}}
{\centerline {78045 Versailles cedex}}
{\centerline {email : \texttt{nperrin@fermat.uvsq.fr}}}

\section*{Introduction} 
\addcontentsline{toc}{section}{Introduction}

Dans la premi\`ere partie de cet article, on \'etudie le sch\'ema de
Hilbert des courbes rationnelles lisses trac\'ees sur une vari\'et\'e
homog\`ene. Soient $G$ un groupe alg\'ebrique simple et simplement
connexe sur ${\bf C}$, $P$ un parabolique de $G$ et $\a\in A_1(G/P)$
une classe d'\'equivalence rationnelle de $1$-cycle. On suppose $\a$
positive (c'est \`a dire coupant positivement les classes de diviseurs
effectifs de $G/P$), alors J.F. Thomsen [T], B. Kim et
R. Pandharipande [KP] ont montr\'e (en particulier) que le sch\'ema de
Hilbert des courbes lisses rationnelles de classe $\a$ trac\'ees sur
$G/P$ est lisse connexe mais pas n\'ecessairement non vide. Dans cette
partie, nous construisons de mani\`ere effective des courbes lisses
sur les vari\'et\'es homog\`enes. Nous retrouvons ainsi le
r\'esultat de J.F. Thomsen, B. Kim et R. Pandharipande par une
m\'ethode compl\`etement diff\'erente et montrons l'existence de
courbes lisses.

\vs 0.2 cm

\noi
{\bf Th\'eor\`eme 1} : \textit{(\i) Le sch\'ema de Hilbert des courbes
rationnelles lisses trac\'ees sur $G/P$ de classe $\a$ positive est
irr\'eductible et lisse.}

\noi
\textit{(\i\i) Si $\a$ est strictement positive (c'est \`a dire
d'intersection strictement positive avec tous les diviseurs effectifs de
$G/P$), alors il existe une courbe rationnelle lisse de classe $\a$
sur $\gsp$ sauf pour $\pu$, $\p^2$ et $\pu\times\pu$.}

\vs 0.2 cm

On donnera \'egalement les cas pour lesquels le sch\'ema de Hilbert des
courbes lisses est non vide lorsque $\a$ est positive non strictement
positive.

Pour d\'emontrer ce th\'eor\`eme on \'etudie les orbites de $G/P$
sous un parabolique $P'$. Ces orbites que nous appelerons $P'$-orbites
g\'en\'eralisent les cellules de Schubert classiques (d\'efinies par
l'action d'un Borel). Elles n'apparaissent pas de fa\c con
syst\'ematique dans la litt\'erature. Seuls quelques cas particuliers
sont d\'ecrits par exemple dans [BGG], [K1], [K2], [T] et [LMS]. Les
m\'ethodes de J.F. Thomsen [T] et B. Kim et R. Pandharipande [KP] pour
prouver la premi\`ere partie du th\'eor\`eme $1$ sont totalement
diff\'erentes de la notre. Ils montrent que le sch\'ema des
applications stables est connexe en reliant toute application \`a une
application du bord par action d'un Borel ou d'un tore maximal, puis
irr\'eductible par lissit\'e. La m\'ethode pr\'esent\'ee ici est plus
directe (pas d'\'etude du bord du sch\'ema des applications stables)
et elle permet de montrer l'existence de morphismes non constants de
$\pu$ dans $G/P$ de classe $\a$ positive et m\^eme de courbes lisses
ce qui n'\'etait pas le cas des m\'ethodes de [T] et [KP]. Pour les
courbes de genre $g>0$, on n'a pour le moment que des r\'esultats
partiels. Citons E. Ballico [B] qui a montr\'e l'irr\'eductibilit\'e
du sch\'ema de Hilbert des courbes de degr\'e $d$ et de genre $g$
trac\'ees sur une quadrique de $\p^n$ pour $n\geq 7$ et $d\geq
2g-1$. Notre m\'ethode permet de montrer, dans le cas o\`u $G$ est
$SL_n$ ou $SO_{2n}$ et si $P$ est un parabolique maximal, que le
sch\'ema des morphismes de degr\'e $d$ d'une courbe de genre $g$ vers
$G/P$ est irr\'eductible d\`es que $d$ est grand devant $g$ (on en
d\'eduit le r\'esultat pour le sch\'ema de Hilbert associ\'e). On peut
\'egalement montrer ce r\'esultat si $G$ est $Sp_{2n}$ et si $P$ est
un parabolique maximal qui ne correspond pas aux sous espaces
totalement isotropes maximaux. Enfin, la m\'ethode de E. Ballico se
g\'en\'eralise aux quadriques lisses de rang $n$ pour $3\leq n\leq
7$. Dans tous ces cas except\'e pour $\p^2$ et $\pu\times\pu$, on sait
\'egalement montrer l'existence de courbes lisses.

\vs 0.1 cm

De plus, les $P'$-orbite que l'on introduit dans cette construction
nous permettent, dans la seconde partie de cet article, de donner une
d\'esingularisation $\pi$ des vari\'et\'es de Schubert \textit{plus
fine} que celle de M. Demazure [D]. En effet, elle sera bijective sur
un plus grand ouvert et une d\'esingularisation de Demazure se
factorise par celle-ci. Cependant, je ne sais pas si elle est
\textit{la plus fine possible} c'est \`a dire si elle est bijective
sur le lieu lisse des vari\'et\'es de Schubert. On peut n\'eanmoins
v\'erifier que c'est le cas pour les groupes de petite dimension
($SL_4$ par exemple). On montre ainsi le r\'esutat suivant :

\vs 0.2 cm

\noi
{\bf Th\'eor\`eme 2} : \textit{(\i) Une d\'esingularisation de Demazure
se factorise par $\pi$.}

\noi
\textit{(\i\i) Le morphisme $\pi$ est une d\'esingularisation des
vari\'et\'es de Schubert.}

\vs 0.2 cm

On donnera aussi une condition suffisante (non
n\'ecessaire) de lissit\'e des vari\'et\'es de Schubert et un
crit\`ere pour qu'une vari\'et\'e de Schubert soit homog\`ene sous
l'action d'un sous groupe de $G$.

\vs 0.5 cm

\centerline{\bf {CONSTRUCTION DE COURBES LISSES SUR LES VARI\'ET\'ES
HOMOG\`ENES}}

\vs 0.5 cm

On adoptera tout au long de cet article les notations de W. Fulton et
J. Harris [FH] pour tout ce qui concerne les groupes, leurs alg\`ebres de
Lie et leurs syst\`emes de racines.
Soit $G$ un groupe de Lie simple et simplement connexe et soit
${\mathfrak{g}}$ son alg\`ebre de Lie. Soit $\gh$ une alg\`ebre de
Cartan de $\gg$ et ${\mathfrak{b}}$ un borel de $\gg$ contenant $\gh$. On note $R$,
respectivement $R^+$ et $R^-$ l'ensemble des racines, respectivement
l'ensemble des racines positives et n\'egatives. On a les
d\'ecompositions de $\gg$ et $\gb$ suivant les poids :
$$\gg=\gh\oplus_{\alpha\in\gh^*}\gg_{\alpha}\ \ \ \  \ \gb=\gh\oplus_{\alpha\in R^+}\gg_{\alpha}$$
Le Groupe de Picard de $G/B$ est exactement le r\'eseau des poids de
$\gh^*$.
Dans le groupe de Picard de $G/B$, le c\^one form\'e par les poids
$x$ tels que : $(x,\alpha)\geq 0$ pour toute racine simple $\a$ (c'est
la chambre de Weyl principale) est engendr\'e par les poids
fondamentaux (diviseurs engendrant le groupe de Picard des $G/P$ o\`u
$P$ est un parabolique maximal). On l'appelle c\^one ample de $\pic(G/B)$.
Le groupe $A_1(G/B)$ est le dual de $\pic(G/B)$. C'est
donc le r\'eseau des racines dans $\gh$. Le c\^one ample de
$\pic(G/B)$ d\'efinit par dualit\'e le c\^one positif de
$A_1(G/B)$. Ce c\^one est le c\^one
engendr\'e par les racines simples. On appelle c\^one strictement
positif les \'el\'ements qui sont strictement positifs sur toutes les
ar\^etes du c\^one ample du groupe de Picard, c'est le c\^one positif
priv\'e de ses faces de codimension $1$.

Un sous groupe parabolique $P$ de $G$ contenant $B$ est donn\'e par un
ensemble de racines simples n\'egatives $\Sigma$ et l'alg\`ebre
correspondante est alors :
$$\gp(\Sigma)=\gb\oplus_{\alpha\in
T(\Sigma)}\gg_{\alpha}$$
o\`u $T(\Sigma)$ est l'ensemble des racines obtenues comme somme de
racines simples en dehors de $\Sigma$. Le groupe de Picard de
$G/P(\Sigma)$ est donn\'e par le sous r\'eseau des poids (c'est \`a
dire le sous r\'eseau de $\gh^*$) d\'efini par les \'equations
$(x,\alpha )=0$ pour toute les racines simples $\alpha$ en dehors de
$\Sigma$. On peut aussi le d\'ecrire comme le r\'eseau engendr\'e par
les coracines ${\check \alpha}$ des racines simples $\alpha\in \Sigma$.
On peut ainsi d\'efinir par restriction un c\^one
ample dans $\pic(G/P)$. Par dualit\'e on peut d\'efinir
des c\^ones positif et strictement positif dans $A_1(G/P)$ qui sont les quotients
des c\^ones de $A_1(G/B)$. 

\vs 0.4 cm

\noi
{\bf Remarque 1} : Si $C$ est une courbe irr\'eductible dans $G/P$ alors
sa classe $[C]$ dans $A_1(G/P)$ est n\'ecessairement dans le c\^one
positif car tous les diviseurs $D$ qui forment une ar\^ete du c\^one ample
de $\pic(G/P)$ sont effectifs et par l'action du groupe on voit que
$(C,D)\geq 0$ (voir par exemple [Kl]). Ainsi [C] est toujours dans le
c\^one positif et la condition de la premi\`ere partie du th\'eor\`eme $1$
est une condition n\'ecessaire pour que le sch\'ema de Hilbert soit
non vide.

\section{Les courbes rationnelles et leur sch\'ema de Hilbert}

On note $\H(\alpha,X)$ le sch\'ema de Hilbert des courbes rationnelles lisses dans $X$ dont dans la classe dans $A_1(X)$ est $\alpha$.

\vs 0.4 cm

\noi
{\bf Proposition 1} : \textit{Soit $\alpha\in A_1(G/P)$ dans le c\^one
positif, alors $\H(\alpha,G/P)$ est lisse de dimension :}
\vs -0.2 cm
$$\sum_{\alpha'\in\Sigma}({\check \alpha'},\alpha )+\ {\rm{dim}}(G/P)-3$$

\vs 0.1 cm

\dm :
Il suffit de montrer que $T_{\gsp}$ est engendr\'e par ses
sections pour avoir la lissit\'e du sch\'ema de Hilbert. Sa dimension au point $C$ est alors donn\'ee par ${\rm{deg}}(T_{\gsp}\vert_C)+\
{\rm{dim}}(\gsp)-3$ c'est ce qu'on cherche car le degr\'e de
$T_{\gsp}\vert_C$ est $\sum_{\alpha'\in\Sigma}({\check \alpha'},\a)$.

Si on a une action lin\'eaire de $P$ sur un espace vectoriel $N$,
alors $P$ agit sur le produit $G\times N$ par
$(p,(g,n))\mapsto(gp^{-1},pn)$. Le quotient $G\times^PN$ est un
fibr\'e vectoriel au dessus de $G/P$. Ainsi, le faisceau $T_{\gsp}$
est le fibr\'e vectoriel associ\'e \`a la repr\'esentation de $P$
suivante : $\gg/\gp$. Mais $G$ agit sur $\gg$ donc le fibr\'e
vectoriel ${\cal G}$ associ\'e \`a $\gg$ est trivial. Comme quotient
du fibr\'e trivial ${\cal G}$, le faisceau $T_{\gsp}$ est engendr\'e
par ses sections. Pour une autre d\'emonstration voir [Ko]
th\'eor\`eme $1.4$ page 241.

\vs 0.4 cm

La difficult\'e de la premi\`ere partie du th\'eor\`eme $1$ r\'eside donc
dans l'irr\'eductibilit\'e du sch\'ema de Hilbert. Pour l'\'etudier,
on va s'int\'eresser au sch\'ema des morphismes $\Mor{\a}{\gsp}$ qui
param\'etrise les morphismes $f$ de sch\'emas de $\pu$ dans $G/P$ tels
que la classe de $f(\pu)$ dans $A_1(G/P)$ est $\a$. L'ensemble des
plongements de $\pu$ dans $\gsp$ de classe $\a$ est un ouvert de ce
sch\'ema. Cet ouvert domine $\H(\a,\gsp)$ (c'est une fibration
en $SL_2$) et l'irr\'eductibilit\'e du sch\'ema de Hilbert se
d\'eduira de celle de $\Mor{\a}{\gsp}$.
La proposition suivante nous permet de ramener ce probl\`eme \`a celui
d'un ouvert :

\vs 0.4 cm

\noi
{\bf Proposition 2} : \textit{Soit $X$ une vari\'et\'e munie d'une
action transitive d'un groupe $G$ et soit $\a\in A_1(X)$. Supposons
qu'il existe un ouvert $U$ de $X$ dont le compl\'ementaire $Z$ est de
codimension sup\'erieure ou \'egale \`a $2$ et tel que $\Mor{\a}{U}$ est
irr\'eductible, alors $\Mor{\a}{X}$ l'est aussi.}

\vs 0.2 cm

\dm :
Les ouverts $\Mor{\a}{g.U}$ de $\Mor{\a}{X}$ recouvrent
$\Mor{\a}{X}$. On utilise pour cel\`a un r\'esultat de Kleiman [Kl] : si $C$
est une courbe de $X$ et si $Z$ est un ferm\'e de codimension au moins
$2$, alors il existe un ouvert de $G$ tel que, pour tout point $g$ de
cet ouvert, l'intersection de $C$ avec les translat\'es $g.Z$ du
ferm\'e $Z$ est vide. Donc si on a un point $f$ dans $\Mor{\a}{X}$ et $C$
son image, alors il existe un ouvert de $G$ tel que $C$ est contenue
dans $g.U$ et donc $f$ est dans $\Mor{\a}{g.U}$ pour tout $g$ dans cet
ouvert. Ceci impose l'irr\'eductibilit\'e. En effet, deux ouverts
$\Mor{\a}{g.U}$ et $\Mor{\a}{g'.U}$ se coupent toujours : il suffit
d'exhiber une courbe qui ne rencontre pas $g.Z\cup g'.Z$. Soit donc
$f\in\Mor{\a}{X}$ quelconque et $C$ son image, le m\^eme r\'esultat de
Kleiman [Kl] nous dit qu'il existe un ouvert de $G$ tel que pour tout
\'el\'ement $g''$ de cet ouvert la courbe $g''.C$ ne rencontre pas
$g.Z\cup g'.Z$. Toutes ces courbes $g''.C$ sont donc dans
l'intersection.
Supposons maintenant que $\Mor{\a}{X}$ a plusieures composantes
irr\'eductibles et soient ${\bf H}$ et ${\bf H'}$ deux telles
composantes. Comme les ouverts $\Mor{\a}{g.U}$ recouvrent et sont irr\'eductibles, il existe deux
\'el\'ements $g$ et $g'$ de $G$ tels que $\Mor{\a}{g.U}$ est un ouvert
non vide de ${\bf H}$ et $\Mor{\a}{g'.U}$ est un ouvert non vide de
${\bf H'}$. Mais alors on sait que $\Mor{\a}{g.U}\cap\Mor{\a}{g'.U}$ est
non vide, c'est donc un ouvert dense de ${\bf H}$ et de ${\bf H'}$ ce
qui impose l'\'egalit\'e de ces composantes.

\vs 0.4 cm

Le principe de la d\'emonstration sera le suivant : on commence par se
rammener grace \`a la proposition pr\'ec\'edente \`a des ouverts de
$\gsp$ dont le compl\'ementaire est de codimension sup\'erieure ou
\'egale \`a $2$ : les $P'$-orbites (voir paragraphe suivant). On
\textit{d\'evisse} ensuite ces $P'$-orbites en vari\'et\'es plus
\textit{simples} pour lesquelles on sait r\'esoudre le
probl\`eme. Deux cas se pr\'esenterons alors selon que $P$ est un
Borel ou non :

- Un fibr\'e en droites projectives au dessus d'une
vari\'et\'e homog\`ene pour laquelle on sait r\'esoudre le
probl\`eme. Dans ce cas on aura besoin d'une condition sur le degr\'e
de la courbe par rapport \`a cette fibration.

- Une tour de fibr\'es affines, de fibr\'es vectoriels
\textit{direction} engendr\'es par leurs sections, au dessus d'un
produit de vari\'et\'es homog\`enes sous des groupes dont le diagramme
de Dynkin est de longueur strictement inf\'erieure \`a celle du
diagramme de Dynkin de $G$.

Pour chacun de ces cas on a une proposition qui permet d'obtenir le
r\'esultat. On commence par d\'efinir ce qu'on appelle une tour de
fibr\'es affines.

\vs 0.4 cm

\noi
{\bf D\'efinition} : Un morphisme $X\stackrel{f}{\fl} Y$ est appel\'e tour de fibr\'es affines si $f$ se d\'ecompose en $X\stackrel{f_1}{\fl} X_1\cdots \stackrel{f_n}{\fl} X_n\stackrel{f_{n+1}}{\fl} Y$ o\`u les $f_i$ sont des fibr\'es affines.

\vs 0.4 cm

\noi
{\bf Proposition 3} : \textit{Si $X\stackrel{\vp}{\fl}Y$ est une tour de
fibr\'es affines de fibr\'es vectoriels direction engendr\'es par leurs
sections, si $\a\in A_1(X)$, alors $\Mor{\a}{X}\fl\Mor{\vp_*\a}{Y}$
est une tour de fibr\'es affines. En particulier, si
$\Mor{\vp_*\a}{Y}$ est irr\'eductible, alors $\Mor{\a}{X}$ l'est
aussi.}

\vs 0.2 cm

\dm :
Il suffit par r\'ecurrence de prouver le cas d'un fibr\'e affine
${\cal F}$ dont le fibr\'e vectoriel direction $F$ est engendr\'e par
ses sections. On a le diagramme suivant :
$$\begin{array}{ccc}
\Mor{\vp_*\a}{Y}\times\pu & \stackrel{p}{\fl} & Y \\
\mapdown{q} & & \\
\Mor{\vp_*\a}{Y} & &
\end{array}$$
Le faisceau $F$ \'etant engendr\'e par ses section le faisceau
$R^1q_*p^*F$ est nul et $q_*p^*F$ est localement libre. 
On va montrer le r\'esultat suivant : le sch\'ema $\Mor{\a}{X}$ au
dessus de $\Mor{\vp_*\a}{Y}$ est le fibr\'e affine de base $q_*p^*F$
associ\'e \`a l'\'el\'ement de $H^1q_*p^*F$ image de celui de
$H^1F$ d\'efinissant ${\cal F}$.

En effet, on a un morphisme universel
$\pu\times\Mor{\vp_*\a}{Y}\fl Y\times\Mor{\vp_*\a}{Y}$ et si on fait
le produit fibr\'e avec $X\times\Mor{\vp_*\a}{Y}$ on obtient une
vari\'et\'e $Z$ qui est le fibr\'e affine $p^*{\cal F}$ au dessus de
$\pu\times\Mor{\vp_*\a}{Y}$. La vari\'et\'e $\Mor{\a}{X}$ est alors
donn\'ee sur chaque fibre au dessus de $\Mor{\vp_*\a}{Y}$ par les
sections de ce fibr\'e, c'est \`a dire par le faisceau d'espaces
affines $q_*p^*{\cal F}$ (d\'efini de la m\^eme fa\c con que pour un
fibr\'e vectoriel). Comme $q_*p^*F$ est localement libre, $q_*p^*{\cal
F}$ est un fibr\'e affine de base $q_*p^*F$, il est associ\'e \`a
l'\'el\'ement de $H^1q_*p^*F=H^1p^*F$ image de celui de $H^1F$
d\'efinissant ${\cal F}$.




\vs 0.4 cm

\noi
{\bf Proposition 4} : \textit{Soit $X\stackrel{\vp}{\fl}Y$ un
fibr\'e en droites projectives de fibr\'e tangent relatif $T_{X/Y}$ et
soit $\a\in A_1(X)$ tel que $\a\cap T_{X/Y}\geq 0$, alors
$\Mor{\a}{X}$ est un ouvert d'un fibr\'e projectif au dessus de
$\Mor{\vp_*\a}{Y}$. En particulier, si $\Mor{\vp_*\a}{Y}$ est
irr\'eductible alors $\Mor{\a}{X}$ l'est aussi.}

\vs 0.2 cm

\dm :
Soit $E$ un faisceau localement libre de rang $2$ sur $Y$ tel que
$X={\rm{Proj}}_Y({\rm{Sym}}(E))$. On a encore un morphisme de
$\Mor{\a}{X}$ dans $\Mor{\vp_*\a}{Y}$. Soit $f\in\Mor{\vp_*\a}{Y}$ et soit
$S={\rm{Proj}}_{\pu}({\rm{Sym}}(f^*E))$ la surface r\'egl\'ee obtenue
comme produit fibr\'e de $\pu$ avec $X$ au dessus de $Y$. On note
$\oo_S(0,1)$ le quotient tautologique de cette surface r\'egl\'ee et
$\oo_S(1,0)$ le diviseur d'une fibre (le groupe de Picard de $S$ est
${\bf Z}^2$). Les rel\`evements de $f$ dans $X$ de classe $\a$
correspondent aux sections de $E(c)$ (o\`u $c$ est un entier
d\'ependant de $\a$) qui sont partout injectives. Les
sections de $E(c)$ correspondent exactement aux sections de
$\oo_S(c,1)$. L'entier $c$ v\'erifie la condition suivante : 
$${\displaystyle{\bigl(c,1\bigr)}}\left(\begin{array}{cc}
0 & 1\\
1 & a\\
\end{array}\right)
\Bigl(\!\!\begin{array}{c}
-a \\
2 \\
\end{array}\!\!\Bigr)=b$$
o\`u la matrice est celle de la forme d'intersection sur $S$, le
vecteur de droite est la classe de $T_{X/Y}$, $a=(\vp_*\a)\cap c_1(E)$ et
$b=\a\cap T_{X/Y}$. On voit ainsi que $c={1\over{2}}(b-a)$
(n\'ecessairement $b\equiv a \ \![2]$). On s'int\'eresse donc au
faisceau $F=q_*(p^*E\otimes\oo_{\pu}(c))$ sur $\Mor{\vp_*\a}{Y}$ (avec
les notations de la proposition pr\'ec\'edente pour l'incidence avec
$Y$). Il est localement libre de rang $b+2$ au dessus de l'ouvert $U$
compl\'ementaire de
${\rm{Supp}}(R^1q_*((p^*E)\otimes\oo_{\pu}(c-1)))$. Cet ouvert
contient l'image de $\Mor{\a}{X}$. En effet, si
$f\in{\rm{Supp}}(R^1q_*((p^*E)\otimes\oo_{\pu}(c-1)))$, alors le
faisceau $f^*E(c)$ a un facteur de degr\'e inf\'erieur \`a $-1$ et ses
sections s'annulent toutes et ne peuvent ainsi donner des courbes sur
$X$. Au dessus d'un point de l'ouvert $U$, les \'el\'ements de la fibre
sont donn\'es par les sections partout injectives qui forment un
ouvert des sections. Sur ${\rm{Proj}}_{U}({\rm{Sym}}({\check
F}))\times\pu$ on a une section tautologique de
$p^*E\otimes\oo_{\pu}(c)$, soit $Z$ la projection dans
${\rm{Proj}}_{U}({\rm{Sym}}({\check F}))$ du lieu des z\'eros de cette
section. On peut maintenant affirmer que $\Mor{\vp_*\a}{Y}$ s'identifie \`a
l'ouvert de ${\rm{Proj}}_{U}({\rm{Sym}}({\check F}))$ compl\'ementaire
$Z$.

Remarquons qu'il est possible que $\Mor{\a}{X}$ soit vide si l'ouvert
image $U$ dans $\Mor{\vp_*\a}{Y}$ est vide.

Cette proposition est encore vraie si on remplace le fibr\'e en
droites projectives par un fibr\'e en espaces projectifs de dimension
sup\'erieure.

\section{D\'ecomposition en $P'$-orbites}

\subsection{D\'efinition et propri\'et\'es g\'en\'erales}

On va ici introduire une classe de vari\'et\'es plus g\'en\'erales que
les cellules de Schubert classiques. Ces vari\'et\'es n'apparaissent
pas ou seulement partiellement dans la litt\'erature : certains
auteurs en ont d\'ecrit des cas particuliers (Kempf [K1] et [K2],
Bernstein, Gelfand et Gelfand [BGG], Lakshmibai, Musili et Seshadri
[LMS], Thomsen [T]). On fixe maintenant un tore $T$ et deux
paraboliques $P$ et $P'$ contenant ce tore et dont l'intersection
contient un Borel (que nous ne fixons pas \`a priori).

\vs 0.4 cm

\noi
{\bf D\'efinition} : On appelle $P'$-orbite de $G/P$ les
orbites de $\gsp$ sous l'action de $P'$.

\vs 0.4 cm

\noi
{\bf Remarques 2} : (\i) Les $P'$-orbites de $\gsp$ forment une
stratification de $\gsp$, elles sont isomorphes \`a $P'wP/P$ ou \`a
$w(P')/(w(P')\cap P)$ o\`u $w\in W$ ($W$ est le groupe de Weyl de
$G$). Si on note pour tout parabolique $Q$, $W(Q)$ le sous groupe du
groupe de Weyl qui laisse stable le parabolique $Q$ (si on fixe un
Borel $B_0$ dans $Q$, ce groupe peut \^etre vu comme le sous groupe de
$W$ engendr\'e par les r\'eflexions par rapport aux racines de $B_0$
dont les oppos\'es sont dans $Q$), alors les $P'$-orbites de $G/P$
sont param\'etr\'ees par les orbites de $W$ sous l'action de
$W(P)\times W(P')$ donn\'ee de la fa\c con suivante :
$(g,(w,w'))\mapsto wg{w'}^{-1}$. La $P'$-orbite $P'wP/P$ est r\'eunion
disjointe de cellules de Schubert, ceci vient de l'\'ecriture 
$$P'wP=\coprod_{(w_1,w_2)\in W(P')\times W(P)}Bw_1ww_2B$$
avec $B$ un Borel de $P\cap P'$.


\noi
(\i\i) Les cellules de Schubert classiques sont d\'ecrites par les
vari\'et\'es $BwB/B$ o\`u $B$ est un Borel de $G$ et $w\in W$. Ce sont
les $B$-orbites de $G/B$.


\noi
Des cas plus g\'en\'eraux ont \'et\'es d\'ecrits par exemple dans
Kempf [K1] et [K2], Berstein, Gelfand et Gelfand [BGG], Lakshmibai,
Musili, Seshadri [LMS] ou Thomsen [T] qui \'etudient les vari\'et\'es
$B'wP/P$ qui sont les cellules de Schubert de $G/P$ par rapport \`a
$B'$ (cas o\`u $P'=B'$ est un Borel contenu dans $P$) ou parfois leurs
\textit{sym\'etriques} : $P'wB/B$ qui sont des $P'$-orbites
plus grandes que les cellules classiques (cas o\`u $P=B$ est un
Borel contenu dans $P'$).

\vs 0.4 cm

\noi
{\bf Exemple 1} : Les $P'$-orbites de $G/P$ sont des ouverts lisses de
$G/P$. Le lemme suivant nous permet de montrer que pour toute
$P'$-orbite $P'\overline{w}P/P$ de $G/P$ ($\overline{w}$ est une
orbite de $W$ sous l'action de $W(P)\times W(P')$), il existe un Borel
$B$ de $P$ et un \'el\'ement $w'\in\overline{w}$ tel que
$P'\overline{w}P/P$ contient la cellule de Schubert $Bw'P/P$ et est
contenue dans l'adh\'erence de cette cellule. On voit ainsi que
l'ad\'erence de la $P'$-orbite $P'\overline{w}P/P$ est une vari\'et\'e
de Schubert.


\vs 0.4 cm

\noi
{\bf Lemme 1} : \textit{Soient $P$ et $P'$ deux paraboliques contenant
un m\^eme Borel, soit $\overline{w}\in W/(W(P)\times W(P'))$, il
existe un Borel $B\subset P$ et un \'el\'ement $w'\in\overline{w}$
tels que $P'\overline{w}P\subset\overline{Bw'B}$.}

\vs 0.2 cm

\dm :
On commence par montrer que l'espace tangent de $P'\overline{w}P$ en
$w$ est $\gp+\gp'$ o\`u $\gp$ est l'espace tangent de $P$ en $e$
(l'identit\'e) et $\gp'$ est l'espace tangent de $w(P')$ en $e$. En
effet, l'espace tangent de $P'we=w.w(P')$ en $w$ est le translat\'e
\`a gauche par $w$ de $\gp'$ et l'espace tangent de $ewP$ en $w$ est
le translat\'e \`a gauche par $w$ de $\gp$. Ainsi $\gp+\gp'$ est
contenu dans l'espace tangent de $P'wP$ en $w$. Mais $P'wP/P$
s'identifie \`a $w(P')/(w(P')\cap P)$ donc $P'wP$ est de dimension
${\rm{dim}}(\gp+\gp')$. La lissit\'e de $P'wP$ nous permet de conclure
que ${\rm{T}}_e(P'wP)=\gp+\gp'$.

Ainsi, notre lemme se ram\`ene \`a montrer qu'il existe deux Borels
$B\subset P$ et $B'\subset w(P')$ tels que $\gb+\gb'=\gp+\gp'$. En
effet, le Borel $w(B)$ est contenu dans $w(P')$ donc il existe $w''\in
W(P')$ tel que $w''w(B)=B'$ donc $B'=w'(B)$ pour
$w'\in\overline{w}$. On aura alors $Bw'B=w'B'B\subset P'wP$ mais comme
ils ont le m\^eme espace tangent, on aura aussi
$P'wP\subset\overline{Bw'B}$. 

On utilise l'abus de notation suivant : si $\a$ est une racine et
$\gp$ l'alg\`ebre de Lie d'un parabolique, on dit que $\a\in\gp$ si
$\a$ est une valeur propre pour l'action du tore $\gh$ sur $\gp$ ou
encore si $\gg_{\a}$ apparait dans la d\'ecomposition de $\gp$. Soient
$\gb$ et $\gb'$ deux Borels de $\gp$ et $\gp'$, on modifie ces Borels
pour obtenir la propri\'et\'e recherch\'ee. Si il existe
$\a\in\gp+\gp'$ (disons dans $\gp$) telle que $\a\not\in\gb+\gb'$
alors il existe un racine simple $\a_0$ de $\gb$ telle que
$-\a_0\in\gp$ et $-\a_0\not\in\gb+\gb'$. En effet sinon pour toute
racine simple $\a_i$ de $\gb$ telle que $-\a_i\in\gp$ on a
$-\a_i\in\gb+\gb'$ i.e. $-\a_i\in\gb'$. Mais alors $\a$ s'\'ecrit
$\a=\sum(-\a_i)$ o\`u les $\a_i$ sont des racines simples de $\gb$
telles que $-\a_i\in\gb'$ donc $\a\in\gb'$ ce qui est absurde.

On remplace maintenant $\gb$ par $s_{\a_0}(\gb)$ et on a : 
$$s_{\a_0}(\gb)+\gb'=(\gb+\gb')\cup\{-\a_0\}$$
on se ram\`ene ainsi \`a des Borels $\gb$ et $\gb'$ tels que
$\gb+\gb'=\gp+\gp'$.

\vs 0.4 cm

Ce lemme montre \'egalement l'existence d'un Borel $B$ et d'un
\'el\'ement $w'\in\overline{w}$ tel que la fl\`eche
$\overline{Bw'B/B}\fl\overline{P'wP/P}$ entre vari\'et\'es de Schubert
est une fibration en $P/B$.

Supposons que $P=B$ est un Borel, les cellules de Schubert classiques
de $G/B$ sont munies d'un ordre partiel (ordre de Bruhat). Les
$P'$-orbites de $G/B$ d\'ecrites par $P'\overline{w}B/B$ avec
$\overline{w}\in W/W(P')$ contiennent comme ouvert dense une cellule
de Schubert classique $Bw'B/B$ qui est donn\'ee par $w'$ (celui du
lemme $1$) qui est l'\'el\'ement maximal de $\overline{w}$ dans
l'ordre de Bruhat. Cet \'el\'ement est unique car les cellules
associ\'ees \`a deux tels \'el\'ements ont la m\^eme adh\'erence :
celle de $P'\overline{w}B/B$.

\vs 0.4 cm

On va utiliser ces $P'$-orbites pour construire des ouverts
de $G/P$ dont le compl\'ementaire est de codimension sup\'erieure ou
\'egale \`a $2$. On commence par montrer un r\'esultat sur ces
$P'$-orbites qui au vu de la proposition $3$ nous permettra de
faire fonctionner la r\'ecurrence.

\vs 0.4 cm

\noi
{\bf D\'efinition} : Les composantes connexes du diagramme de Dynkin de
$G$ priv\'e de $P'$ sont les diagrammes de Dynkin obtenus \`a partir
de celui de $G$ en enlevant les sommets correspondant \`a $P'$.

\vs 0.4 cm

\noi
{\bf Exemple 2} : Dans $SL_4$, si $P'$ est le parabolique
correspondant aux droites, il y a deux composantes connexes du
diagramme de Dynkin de $G$ priv\'e de $P'$. Ces deux composantes sont
isomorphes au diagramme de Dynkin de $SL_2$.

\vs 0.4 cm

\noi
{\bf Remarques 3} : (\i) Si $N$ est un espace vectoriel muni d'une
action de $P$, alors on peut d\'efinir un fibr\'e vectoriel ${\cal N}$
sur $\gsp$ \`a partir de $N$. Ce fibr\'e est le produit contract\'e
$G\times^PN$. De la m\^eme fa\c con, si $X$ est une vari\'et\'e affine
sur laquelle $P$ agit \`a gaughe, on peut d\'efinir le produit
contract\'e $G\times^PX$ qui est le quotient de $G\times X$ par $P$
(voir par exemple [DG]). Si l'action de $P$ sur $X$ se prolonge en une
action de $G$, alors $G\times^PX$ (et en particulier ${\cal N}$ dans
le cas d'un fibr\'e) est trivial sur $\gsp$.

\noi
(\i\i) Soit $w\in W$. On note $N'$ la partie unipotente de $w(P')$ et
$R'$ le quotient r\'eductif correspondant. La d\'ecomposition de Levi
(voir [Bo] 11.22 et 14.17-19) nous permet de dire que $w(P')$ est le
produit semi-direct de $R'$ et $N'$. Cette situation est rigidifi\'ee
par le choix du tore maximal $w(T)$ contenu dans $w(P')$ : la
d\'ecomposition de $\gg$ en sous espaces propes induit une section $s$
de $R'$ dans $w(P')$. On note $N=N'\cap P$ et $R=R'\cap P$ (c'est
l'image de $P$ dans $R'$). Le groupe $P\cap w(P')$ est produit
semi-direct de $R$ et $N$. On a ainsi deux fibr\'es vectoriels ${\cal
N'}$ et ${\cal N}$ sur $R'/R$.

\noi
(\i\i\i) Le groupe $R'$ est produit des groupes $G_i$ donn\'es par les
composantes du diagramme de Dynkin de $G$ priv\'e de $P'$  et de
facteurs isomorphes \`a ${\bf C}^*$ (il y en a
${\rm{rang}}({\rm{Pic}}(G/P'))$). De m\^eme, $R$ est produit de
paraboliques $P_i$ des $G_i$ et des m\^emes facteurs isomorphes \`a
${\bf C}^*$. Les $P_i$ sont donn\'es par les sommets de $P$ regard\'es
dans le diagramme de Dynkin de $G$ prov\'e de $P'$.

\vs 0.4 cm


\noi
{\bf Proposition 5} : \textit{Le morphisme naturel $f$ de $P'wP/P$
vers $R'/R$ est une tour de fibr\'es affines dont les fibr\'es
vectoriels direction correspondent \`a la d\'ecomposition naturelle de
${\cal N'}/{\cal N}$. Ils sont d\'efinis sur $R'/R$ et sont
engendr\'es par leurs sections.}

\vs 0.2 cm

\dm :
On regarde ici la $P'$-orbite $P'wP/P$ sous la forme $w(P')/(w(P')\cap
P)$. 
Soit $1=N'_0\subset N'_1\subset\cdots\subset N'_n=N'$ la suite centrale
ascendante de l'unipotent $N'$. On note $Z'_i=N'_i/N'_{i-1}$ qui est le
centre de $N'/N'_{i-1}$. C'est un espace vectoriel sur lequel $R'$
agit lin\'eairement. Posons $P'_i=w(P')/N'_i$, $P_i$ l'image de $P$
dans $P'_i$ et $Z_i=P_i\cap Z'_i$ qui est un espace vectoriel sur
lequel $R$ agit lin\'eairement. L'\'ecriture de $f$ comme tour de
fibr\'e affine est 
$w(P')/(w(P')\cap P)=P'_0/P_0\fl P'_1/P_1\fl\cdots
\fl P'_i/P_i\fl\cdots\fl P'_n/P_n=R'/R$.
Les fl\`eches sont des fibrations
localement triviales pour la topologie de Zariski et la fibre de
$P'_{i-1}/P_{i-1}\fl P'_i/P_i$ au dessus du point $1$ est
$Z'_i/Z_i$. Il sagit de montrer que cette fibration est affine de
fibr\'e vectoriel direction l'image r\'eciproque (par la projection
$P'_i/P_i\fl R'/R$) du fibr\'e vectoriel sur $R'/R$ d\'eduit de
la repr\'esentation $Z'_i/Z_i$ de $R$. 

Notons $A=P'_{i-i}$, $B=P_{i-1}$ et $U=Z'_i$. La fibration
$P'_{i-1}/P_{i-1}\fl P'_i/P_i$ est donc $A/B\fl (A/U)/(B/(U\cap
B))=A/BU$. On va montrer le lemme suivant :

\vs 0.2 cm

\noi
{\bf Lemme 2} : \textit{Soient $A$ un groupe lin\'eaire, $B$ et $U$
des sous groupes ferm\'es de $A$. On suppose $U$ unipotent commutatif
dans $A$, alors la projection $A/B\stackrel{p}{\fl} A/BU$ est un
fibr\'e affine dont le fibr\'e direction sur $A/BU$ est d\'eduit de la
repr\'esentation (par automorphismes int\'erieurs) de $BU$ sur
$U/(U\cap B)$.}

\vs 0.2 cm

\dm :
On pose $C=BU$ et on consid\`ere la vari\'et\'e $A\times (C/B)$ qui est
munie d'une action de $C$ donn\'ee de la fa\c con suivante :
$(c,(a,x))\mapsto (ac^{-1},cx)$. On peut alors regarder le quotient
$A\times^C (C/B)$ qui est muni d'une projection vers $A/C$. La fibration
ainsi obtenue : $A\times^C (C/B)\fl A/C$ est exactement la fibration $p$
de l'\'enonc\'e. En effet, si on regarde l'action restreinte de $B$
sur $A\times (C/B)$, le quotient $A\times^B (C/B)$ a une projection vers
$A/B$ qui est scind\'ee (il suffit de prendre la section qui \`a
$\widehat{a}$ la classe de $a\in A$ dans $A/B$ associe
la classe de $(a,\overline{1})$ dans $A\times^B (C/B)$, o\`u on a not\'e
$\overline{1}$ la classe de l'identit\'e dans $C/B$). Mais alors cette
section et les morphismes naturels $A\times^B (C/B)\fl A\times^C (C/B)\fl
A/C$ nous donnent la fl\`eche de $A/B$ dans $A/C$ qui a $\widehat{a}$
associe $\tilde{a}$ pour tout $a\in A$. Ainsi, on a une fl\`eche qui
respecte les morphismes vers $A/C$ de $A/B$ dans $A\times^C (C/B)$ et on
peut construire sa r\'eciproque : \`a la classe de $(a,\overline{c})$
dans $A\times^C (C/B)$ on associe $\widehat{ac}$. 

Ainsi $p$ est une fibration de groupe structural $C$ agissant sur
$C/B=U/(B\cap U)$. L'action de $C$ est l'action
naturelle de $C$ sur $C/B$. Mais si $b\in B$ et $u\in U$ alors la
classe de $bu\in C$ dans $U/(B\cap U)$ est celle de $bub^{-1}$ donc
l'action de $C$ sur $U/(B\cap U)$ est donn\'ee par
$(bu,\overline{u'})\mapsto \overline{bub^{-1}.bu'b^{-1}}$. On voit ainsi que
cette action est affine : la partie $bu'b^{-1}$ \'etant lin\'eaire et
on a la translation par $bub^{-1}$. La fibration $p$ est donc affine
et sa direction est donn\'ee par la partie lin\'eaire de l'action
c'est \`a dire par l'action de $C$ sur $U/(B\cap U)$ donn\'ee par
$(bu,\overline{u'})\mapsto \overline{bu'b^{-1}}$. C'est l'action
restreinte de $B\subset C$ sur $U/(B\cap U)$ par automorphisme
int\'erieur.

\vs 0.2 cm

Il nous reste \`a voir que tous les fibr\'es vectoriels associ\'es \`a
ces fibr\'es affines sont d\'efinis sur $R/R'$ et qu'ils sont tous
engendr\'es par leurs sections. Il suffit de montrer que dans le cas
de la fibration $P'_{i-1}/P_{i-1}\fl P'_i/P_i$ l'action de $P_i$
sur $Z'_i/Z_i$ est en fait donn\'ee par une action de $P_{i+1}$. Mais
l'action est donn\'ee par $(p,\overline{z'})\mapsto
\overline{pz'p^{-1}}$. Si on remplace $p$ par $pz$ avec $z\in Z_{i+1}$
alors on a $\overline{pzz'z^{-1}p^{-1}}=\overline{pz'p^{-1}}$ car $N$
est ab\'elien. Si on prend un \'el\'ement $\widetilde{p}$ de
$P_{i+1}$, on peut d\'efinir son action sur $Z'_i/Z_i$ par l'action de
$p$. On voit ainsi que tous les fibr\'es direction sont d\'efinis sur la base
$R'/R$. Enfin, comme l'action de $R$ sur $Z'_i$ se prolonge \`a
$R'$, les fibr\'es associ\'es \`a la repr\'esentation $Z'_i/Z_i$ de
$R$ sont quotient du fibr\'e trivial associ\'e \`a la repr\'esentation
de $R$ (qui se prolonge \`a $R'$) dans $Z'_i$. Ils sont
donc engendr\'es par leurs sections.

Si $N'$ est ab\'elien, le morphisme $P'wP/P\fl R'/R$ est un
fibr\'e vectoriel. En effet, la tour se r\'eduit \`a un fibr\'e affine
et de plus la section de $R'$ dans $P'$ nous d\'efinit une section de
$R'/R$ dans $P'wP/P$ qui nous dit que ce fibr\'e affine est
vectoriel.


%

\subsection{Etude de la codimension du bord de la $P'$-orbite
maximale}

On appelle $P'$-orbite maximale de $G/P$ la $P'$-orbite
$P'wP/P$ qui est dense dans $G/P$. Cette orbite est
unique. On cherche \`a quelle condition cette orbite maximale a un
compl\'ementaire de codimension sup\'erieure ou \'egale \`a $2$. On
commence par donner une condition pour que $P'wP/P$ soit
la $P'$-orbite  maximale de $G/P$. On fixe un Borel $B$ dans $P\cap
P'$.

\vs 0.4 cm

\noi
{\bf Lemme 3} : \textit{La $P'$-orbite $P'wP/P$ est
dense dans $G/P$ si et seulement si $\gp\cap (-w(\gp'))$ contient
l'alg\`ebre de Lie d'un Borel, c'est \`a dire, si et seulement si $w$
est dans l'orbite de $w_0$ (l'\'el\'ement de longueur maximale) sous
$W(P)\times W(P')$.} 

\vs 0.2 cm

\dm : 
Il suffit de d\'emontrer que cette condition implique que l'orbite est
maximale pour l'inclusion car par unicit\'e de cette orbite toutes
les orbites ainsi obtenues seront isomorphes. Les $P'$-orbites sont
d\'ecrites par $w(P')/(w(P')\cap P)$ pour $w\in W$. Cette $P'$-orbite
est isomorphe \`a $P'/(P'\cap w^{-1}(P))$. Pour maximiser cette orbite
on cherche \`a minimiser $P'\cap w^{-1}(P)$. Pour que cette intersection
soit minimale il faut et il suffit que $w^{-1}(\gp)$ contienne le moins de
racines de $\gp'$ possible. C'est le cas si et seulement si $w(\gp')$
contient toutes les racines qui ne sont pas dans $\gp$. Donc si $\gp$
et $-w(\gp')$ contiennent le m\^eme Borel l'orbite est maximale pour
l'inclusion et r\'eciproquement.

\vs 0.4 cm

Soit $B$ un Borel de $G$ et soit $P$ un parabolique de
$G$ contenant $B$, on note $\S(\gp,\gb)$ les racines simples de $\gb$
correspondant aux sommets du diagramme de Dynkin qui d\'efinissent
$\gp$. Par exemple $\S(\gb,\gb)$ est l'ensemble des sommets du
diagramme de Dynkin. On d\'efinit une involution $i$ du diagramme de
Dynkin de la fa\c con suivante : soit $B$ un Borel de $G$ et $\gb$ son
alg\`ebre de Lie, soit $w_0\in W$ le seul \'el\'ement du groupe de
Weyl qui envoie $\gb$ sur $-\gb$ (c'est l'\'el\'ement de longueur
maximale), soit $\a$ une racine simple de $\gb$ (cette racine
correspond \`a un sommet du diagramme de Dynkin), on d\'efinit $i(\a)$
comme \'etant la racine simple de $\gb$ \'egale \`a $-w_0(\a)$. Cette
involution correspond \`a l'involution classique du diagramme de
Dynkin de $A_n$, $D_{2n}$ et $E_6$ et \`a l'identit\'e sur les autres
diagrammes.

\vs 0.4 cm

\noi
{\bf Proposition 6} : \textit{La $P'$-orbite maximale $P'w_0P/P$
a un compl\'ementaire de codimension sup\'erieure ou \'egale \`a $2$
dans $\gsp$ si et seulement si dans le diagramme de Dynkin, $\S(\gp,\gb)$ et
$i(\S(\gp',\gb))$ sont disjoints.}

\vs 0.2 cm

\dm :
Il suffit de montrer que l'application naturelle $p$ de $\pic(\gsp)$ dans
$\pic(P'w_0P/P)$ est injective. Or le noyau de cette
application est donn\'e par  $\pic(\gsp)\cap \pic(G/w_0(P'))$ dans
$\pic(G/(P\cap w_0(P')))\subset \gh^*$ (car le noyau de l'application
$\pic(G/(P\cap w_0(P'))) \fl \pic(w_0(P')/(P\cap w_0(P')))$ est
$\pic(G/w_0(P'))$).

Pour que $p$ soit injective, il faut et il suffit que cette
intersection soit nulle. Or $\pic(G/P)$ est l'orthogonal dans $\gh^*$
de l'ensemble $\a(\gp)$ des racines $\a\in\gp$ telles que $-\a\in\gp$. On
voit que l'intersection $\pic(\gsp)\cap \pic(G/w_0(P'))$ est nulle si et
seulement si $\a(\gp)\cup\a(w_0(\gp'))$ engendre tout $\gh$. Si $RS$
est l'ensemble des racines simples de $\gb$, alors
$\S(\gp,\gb)=RS\setminus(\a(\gp)\cap RS)$ et
$\S(-w_0(\gp'),\gb)=RS\setminus(\a(w_0(\gp'))\cap RS)$. On voit alors que
$\a(\gp)\cup\a(w_0(\gp'))$ engendre tout $\gh$ si et seulement si
$(\a(\gp)\cap RS)\cup(\a(w_0(\gp'))\cap RS)=RS$ ce qui est \'equivalent \`a
$\S(\gp,\gb)$ et $\S(-w_0(\gp'),\gb)=i(\S(\gp',\gb))$ sont disjoints.

\vs 0.4 cm

\noi
{\bf Remarques 4} : (\i) Cette condition nous permet de construire
pour tout parabolique $P$ qui n'est pas un Borel un parabolique $P'$
tel que la $P'$-orbite $P'w_0P/P$ soit maximale et que son
compl\'ementaire soit de codimension au moins $2$. En effet, soit $P$
un parabolique qui n'est pas un Borel, alors il correspond dans le
diagramme de Dynkin \`a un ensemble $\S$ de sommets qui ne contient
pas tous les sommets, il suffit alors de prendre pour $P'$ un parabolique
maximal dont le sommet $s$ dans le diagramme de Dynkin n'est pas dans
$i(\S)$. Par contre cette proposition nous montre que ceci ne sera jamais
possible avec les Borels. On propose donc une autre m\'ethode pour
r\'esoudre ce cas (voir proposition $7$).

\noi
(\i\i) Soit $P$ un parabolique qui n'est pas un Borel et $P'$ un
parabolique tel que $P$ et $P'$ v\'erifient les hypoth\`eses de la
proposition $6$. Alors la vari\'et\'e $R'/R$ obtenue \`a la
proposition $5$ est donn\'ee par le produit des vari\'et\'es
homog\`enes (sous les groupes $G_i$ definis par les composantes
connexes du diagramme de Dynkin de $G$ priv\'e de $P'$) de paraboliques
 $P_i$ d\'efinis par les points du diagramme de Dynkin de
$i(\S(\gp))$. Autrement dit, si on consid\`ere sur le diagramme de
Dynkin les sommets de $P$ apr\`es involution et ceux de $P'$ alors ces
diagrammes sont disjoints et les $P_i$ sont donn\'es
dans les composantes connexes du diagramme de Dynkin de $G$ priv\'e de
$P'$ par les points du diagramme de $i(\S(\gp))$ quand on a retir\'e
ceux de $P'$. Par exemple si on consid\`ere la vari\'et\'e des droites
de $\p^4$ et que l'on prend pour $P'$ un parabolique fixant une droite
alors $R'/R$ est $\p^2$.


\vs 0.4 cm

\noi
{\bf Proposition 7} : \textit{Soit $B$ un Borel et $P$ un parabolique contenant $B$ dont le diagramme de Dynkin (par rapport \`a $B$) a trois sommets
cons\'ecutifs et que celui du milieu n'est rattach\'e qu'\`a deux
sommets dans le diagramme de Dynkin de $G$ (respectivement a deux sommets cons\'ecutifs dont
l'un est au bord du diagramme), soit $P'$ le parabolique obtenu en
enlevant le sommet du milieu (respectivement celui du bord), alors
$G/P$ est un $\pu$-bundle au dessus de $G/P'$}

\vs 0.2 cm

\dm :
Soient $\gp$ et $\gp'$ les alg\`ebres de Lie de $P$ et $P'$, on voit
que $\gp$ est une sous alg\`ebre de $\gp'$ et que si $\a$ est la
racine simple de $B$ qui correspond au sommet que l'on a retir\'e
alors $\gp'$ est l'alg\`ebre de Lie engendr\'ee par $\gp$ et
$-\a$. Mais comme le sommet correpondant \`a $\a$ forme une composante
connexe du diagramme de Dynkin de $G$ priv\'e de $P'$, alors on voit
que $\a$ est orthogonale \`a toutes les racines simples n\'egatives de
$\gp$ ce qui impose $\gp'=\gp\oplus\gg_{\a}$.

\vs 0.4 cm

\noi
{\bf Remarque 5} : Soient $P$ et $P'$ comme dans la proposition
pr\'ec\'edente, si $C$ est une courbe trac\'ee sur $G/P$ dont la
classe est $x\in A_1(G/P)$, alors son degr\'e par rapport \`a la
fibration est donn\'e par $(x,\a)$ o\`u $\a$ est la racine simple
(et donc le caract\`ere) qui correspond au sommet que l'on a
retir\'e. En d'autres termes, avec les notations de la proposition
pr\'ec\'edente, si $\a$ est le sommet (la racine) que l'on a retir\'e,
alors le fibr\'e tangent relatif $T_{G/P/G/P'}$ est de premi\`ere
classe de Chern $\a$ (voir par exemple [D]).


Si $P$ est un Borel, on peut appliquer cette proposition \`a tous les
sommets du diagramme de Dynkin. Si de plus $C$ est une courbe
trac\'ee dans $G/B$ dont la classe $x\in A_1(G/B)$ est dans le c\^one
positif, alors il existe au moins un sommet du diagramme pour lequel
le degr\'e de $C$ sera positif. En effet, $x$ est dans le c\^one
positif si et seulement si pour toute racine simple $\a$ on a
$(x,{\check \a})\geq 0$ (o\`u ${\check \a}$ est la coracine de
$\a$). Mais les coracines forment le c\^one ample et sont donc
combinaisons lin\'eaires \`a coefficients positifs des racines
simples. Si pour toute racine simple $\a$ on a $(x,\a)<0$ alors
$x$ ne peut pas \^etre dans le c\^one positif. Ceci nous permettra donc
d'appliquer la proposition $4$. De m\^eme, si $x$ est dans le c\^one
strictement positif, il existe un sommet du diagramme de Dynkin tel
que $(x,\a)>0$.

\vs 0.4 cm

{\bf Application : d\'emonstration de l'irr\'eductibilit\'e du
sch\'ema de Hilbert} : On raisonne par r\'ecurrence sur la longueur du diagramme de Dynkin (nombre de
sommets). Le cas de $SL_2$ est \'evident. Soit $P$ un
parabolique de $G$ et soit $\a\in A_1(G/P)$ dans le c\^one positif. Si
$P$ est un Borel, on a vu \`a la remarque $6$ qu'il existe un sommet
du diagramme de Dynkin correspondant \`a la racine simple $\a'$ tel que
$(\a,\a')\geq 0$. Dans ce cas $G/B$ est une fibration en droites
projectives au dessus de
$G/P'$ tel que le degr\'e de $\a$ par
rapport \`a cette fibration est positif. La proposition $4$ nous
permet donc de nous ramener au cas o\`u $P$ n'est pas un Borel. 

La remarque $5$ nous permet alors de construire un parabolique $P'$
tel que la $P'$-orbite $P'w_0P/P$ est maximale et que son
compl\'ementaire est de codimension au moins $2$. La proposition $2$
nous permet de nous ramener au probl\`eme d'irr\'eductibilit\'e du
sch\'ema des morphismes pour cette orbite.

Enfin, la proposition $5$ nous dit que $P'w_0P/P$ est
une tour de fibr\'es affines associ\'es \`a des fibr\'es vectoriels
engendr\'es par leurs sections au dessus d'un produit de vari\'et\'es
homog\`enes sous des groupes dont le diagramme de Dynkin est de longueur
strictement inf\'erieure \`a celle du diagramme de Dynkin de $G$. On
conclue par hypoth\`ese de r\'ecurrence en utilisant la proposition
$3$.

\vs 0.4 cm

\noi
{\bf Remarque 6} : On sait de la m\^eme fa\c con que tous les
sch\'emas $\Mor{\a}{U}$ o\`u $U$ est une $P'$-orbite de $G/P$ sont
irr\'eductibles.

\subsection{Existence de courbes lisses}

On va dans ce paragraphe montrer l'existence de courbes lisses sur les
vari\'et\'es homog\`enes. On se restreindra au cas ou la classe de la
courbe est dans le c\^one strictement positif et on
verra (remarque $8$) comment se ramener \`a ce cas si la courbe est
dans le c\^one positif. On appelle courbe nodale une courbe
irr\'eductible et g\'en\'eriquement r\'eduite qui est lisse ou qui a
pour seules singularit\'es des points doubles ordinaires. Si $F$ est un
fibr\'e vectoriel sur un sch\'ema $X$, on dit que $F$ s\'epare les
points si pour tout couple de points $(P,Q)$ de $X$, on a
$h^0(F\otimes\I_P) >h^0(F\otimes\I_{P\cup Q})$. Si $X$ est une
vari\'et\'e et $\a\in A_1(X)$, on note $\mathfrak{H}_{0,\a}(X)$ le
sch\'ema de Hilbert des courbes rationnelles nodales de classe
$\a$. On dit que $X$ v\'erifie la condition $(*)$ si la vari\'et\'e
d'incidence $\{(x,C)\in
X\times\mathfrak{H}_{0,\a}/x\in{\rm{Sing}}(C)\}$ est irr\'eductible.

\vs 0.4 cm

\noi
{\bf Lemme 4} : \textit{Soit $X$ une vari\'et\'e munie d'un fibr\'e
affine ${\cal F}$, de fibr\'e vectoriel direction $F$ et d\'efini par
$\eta\in H^1(X,F)$. Soit $C$ une courbe rationnelle irr\'eductible et
g\'en\'eriquement r\'eduite sur $X$ et $f:\pu\fl C$ une
d\'esingularisation de $C$. On note $\overline{\eta}$ la restriction
de $\eta$ \`a $H^1(C,F\vert_C)$. Si on suppose que $F$ est engendr\'e par ses
sections, alors il existe un rel\`evement de $f$ dans
${\rm{Aff}}({\cal F}\vert_C)$ le fibr\'e affine associ\'e \`a ${\cal
F}$ au dessus de $C$. Si on suppose de plus que $C$ est nodale et que
l'une des conditions suivante est v\'erifi\'ee :}

\textit{(\i) $\overline{\eta}\not =0$ et $X$ v\'erifie $(*)$}

\textit{(\i\i) $f^*F\vert_C$ s\'epare les points}

\noi
\textit{alors il existe un rel\`evement lisse de $f$ dans ${\rm{Aff}}({\cal
F})\vert_C$.}

%
%
%


\vs 0.2 cm

\dm : 
Comme $F$ est engendr\'e par ses sections, $H^1f^*F\vert_C$ est nul donc
$f^*{\cal F}\vert_C$ est le fibr\'e vectoriel $f^*F\vert_C$ qui est
engendr\'e par ses sections. Une telle section nous donne un
rel\`evement de $f$ dans ${\rm{Aff}}({\cal F}\vert_C)$.
Si $C$ est lisse, alors un rel\`evement quelconque de $C$ est lisse.

Si $\overline{\eta}\not =0$, le fibr\'e affine n'est pas vectoriel et n'a donc
pas de section. Prenons un rel\`evement $f'$ de $f$ donn\'e par une
section de $f^*F\vert_C$. Ce rel\`evement est n\'ecessairement non bijectif
sur $C$ (sinon ce serait une section de ${\cal F}\vert_C$). C'est donc une
d\'esingularisation partielle de $C$. Ainsi, il existe un rel\`evement
$f'$ de $f$ qui d\'esingularise au moins un point de $C$. Mais par
monodromie (le groupe de monodromie agit transitivement sur les points
singulier grace \`a la condition $(*)$, cf. [ACGH] ou [Har]) on sait
que pour chaque point double il existe une section de $f^*F\vert_C$ qui
d\'esingularise ce point. En prenant une section g\'en\'erale on
obtient une d\'esingularisation en tout point.

Supposons que $f^*F\vert_C$ s\'epare les points. On sait que $H^1f^*F\vert_C$ est
nul c'est \`a dire $f^*{\cal F}\vert_C$ est le fibr\'e vectoriel
$f^*F\vert_C$. Soit $P$ un point double de $C$ et $x$ et $y$ ses
ant\'ec\'edents par $f$. Il existe une section de $f^*F\vert_C$ qui
s\'epare $x$ et $y$. Mais alors cette section nous donne un
rel\`evement $f'$ de $f$ qui d\'esingularise $P$. Ainsi, il existe un
rel\`evement d\'esingularisant chaque point double et une section
g\'en\'erale de $f^*F\vert_C$ nous donne un rel\`evement lisse.

Si $\overline{\eta}=0$ et qu'on ne suppose plus que $f^*F\vert_C$
s\'epare les points alors il n'existe pas n\'ecessairement de
rel\`evement lisse de $f$. C'est le cas si $F$ est trivial sur une
courbe nodale.

\vs 0.4 cm

\noi
{\bf Remarques 7} : (\i) Avec les notations du lemme, supposons que
$C$ est contenue dans une vari\'et\'e homog\`ene $X$ et que sa classe
dans $A_1(X)$ est dans le c\^one strictement positif. Supposons de
plus que $c_1(F)$ est non nul dans le c\^one ample de $\pic (X)$ alors
il existe un rel\`evement lisse de $f$ dans ${\rm{Aff}}({\cal
F})$. En effet, il suffit de v\'erifier que l'une des condition du
lemme est v\'erifi\'ee. Il suffit donc de v\'erifier que $f^*F$
s\'epare les points. Mais le degr\'e de $f^*F$ sur $\pu$ est
strictement positif donc $f^*F$ s\'epare les points.

\noi
(\i\i) On aura dans la suite besoin de savoir que $\p^2$ v\'erifie la
condition $(*)$. Ceci est fait dans [ACGH].



\vs 0.4 cm 

\noi
{\bf Lemme 5} : \textit{Soient $C$ une courbe rationnelle
irr\'eductible et g\'en\'eriquement r\'eduite et
$\p_C(E)\stackrel{\vp}{\fl} C$ une fibration en droites projectives au
dessus de $C$ d'espace tangent relatif $T$. Soit $f:\pu\fl C$ une
d\'esingularisation de $C$. Supposons que
$f^*E=\oo_{\pu}\oplus\opu(x)$ avec $x\geq 0$. Soit $d\geq 0$ un entier tel
que $d\equiv x\ \![2]$ et $d\geq x$, alors il existe un rel\`evement
$f'$ de $f$ tel que $f'_*[\pu]\cap T=d$. Supposons de plus
que $C$ est nodale et $d>0$, alors on peut choisir $f'$ d'image
lisse.}

\vs 0.2 cm

\dm :
Les rel\`evements de degr\'e relatif $d$ de $f$ sont donn\'es par les
quotients isomorphes \`a $\opu({x+d\over{2}})$ de $f^*E$.Un tel
quotient existe si et seulement si $d\equiv x\ \![2]$ et
${x+d\over2}\geq x$ c'est \`a dire $d\geq x$ ce qui est le cas. En
effet, si les conditions sont v\'erifi\'ees, il existe un tel
quotient. Si $d<x$ alors pour avoir un tel quotient il faut que
${x+d\over2}=0$ c'est \`a dire $d=x=0$. C'est absurde.

Dans le cas o\`u $C$ est nodale , il suffit de v\'erifier que l'on
peut s\'eparer les points. La donn\'ee d'un quotient de $f^*E$
isomorphe \`a $\opu({x+d\over{2}})$ \'etant \'equivalente \`a la
donn\'ee d'une section partout non nulle de $f^*E({d-x\over2})$, on
est ramen\'e \`a montrer que pour tout couple de points $(x,y)$ de $\pu$ il
existe une telle section $s$ telle que $s(x)$ et $s(y)$ sont
lin\'eairement ind\'ependants. Mais
$f^*E({d-x\over2})=\opu({x+d\over{2}})\oplus\opu({d-x\over{2}})$ donc
ceci est possible d\`es que ${d-x\over2}\geq 0$ et
${x+d\over{2}}>0$, ce qui est v\'erifi\'e sous nous hypoth\`eses. 

Si on ne suppose plus $d>0$, il n'existe pas n\'ecessairement de
rel\`evement lisse de $f$. En effet, si $E$ est trivial et $C$ nodale
alors $\p_C(E)=C\times\pu$ n'a pas de rel\`evement lisse de $C$.

La condition $d\equiv x\ \![2]$ ne d\'epend que de la classe $\a$ du
rel\`evement $f'$ et de $c_1(E)$. En effet, on a $d=\a\cap
T$ et $x=(\vp_*\a)\cap c_1(E)$. La condition $d\geq x$ est
\'equivalente \`a $h^1f^*E({d-x\over2}-1)=0$. C'est donc une condition
ouverte.

%

\vs 0.4 cm

{\bf D\'emonstration de l'existence de courbes lisses} : On
proc\`edera de la fa\c con suivante : on commence par supposer que
pour $\a$ positive dans $A_1(G/P)$, le sch\'ema $\Mor{\a}{\gsp}$ est
non vide. On ram\`ene le cas d'un Borel \`a celui d'un parabolique
qui n'est pas un Borel. On montre ensuite le r\'esultat par
r\'ecurrence sur la longueur du diagramme de Dynkin. On initialise la
r\'ecurrence en montrant les cas des groupes dont le diagramme est de
longueur au plus $3$. On montrera le cas de $G_2$ par une m\'ethode
diff\'erente. Enfin on montre que $\Mor{\a}{\gsp}$ est non vide si
$\a$ est positive. 


On suppose connu les r\'esultats suivants : il existe des courbes
rationnelles lisses sur $\p^n$ d\`es que $n\geq 3$, sur
$\pu\times\p^n$ d\`es que $n\geq 2$ (lemme $5$) et sur
$\pu\times\pu\times\pu$ (lemme $5$). Sur $\p^2$ et $\pu\times\pu$ il
existe des courbes nodales. Soit $\a$ dans le c\^one strictement
positif de $A_1(G/P)$.

LE CAS DES BORELS : Si $P=B$ est un Borel, la remarque $5$ nous permet de
trouver un parabolique $P'$ tel que $G/B\fl G/P'$ est une fibration en
droites projectives de fibr\'e tangent relatif $T$ v\'erifiant $\a\cap
T>0$. Comme on a suppos\'e $\Mor{\a}{G/B}$ non vide, il existe une
courbe de $G/P'$ telle que la fibration v\'erifie les conditions du
lemme $5$. Comme ces conditions sont ouvertes, la courbe g\'en\'erale
v\'erifie ces conditions. On est donc ramen\'e \`a prouver l'exitence
de courbes lisses (et m\^eme seulement \`a points doubles ordinaires)
sur $G/P'$. Ceci nous permet notamment de dire que sur $SL_3/B$ il
existe des courbes lisses dont la classe est quelconque dans le c\^one
strictement positif.

LE CAS G\'EN\'ERAL : On proc\`ede par r\'ecurrence sur la longueur du
diagramme de Dynkin. On choisit un parabolique $P'$ tel que la
$P'$-orbite $P'w_0P/P$ est maximale. La proposition $5$ et le
lemme $4$ nous permettent de construire des courbes lisses sur les
ouverts $P'w_0P/P$ : d\`es qu'il y a des courbes lisses sur
$R'/R$, il y en a sur  la $P'$-orbite (on prend une section du fibr\'e
affine). Par hypoth\`ese de r\'ecurrence, les
seuls cas qui posent probl\`eme sont donc ceux o\`u $R'/R$ est $\pu$,
$\p^2$ ou $\pu\times\pu$. Dans ce cas le rang du groupe de Picard
est au plus $2$ et $P$ est donn\'e par au plus deux points
dans le diagramme de Dynkin. On choisit pour
$P'$ le parabolique maximal qui correspond \`a un point extr\^eme du
diagramme de Dynkin. Ainsi dans le cas de $SL_n$ on se ram\`ene \`a
$SL_{n-1}$ ; dans le cas de $Sp_{2n}$ on se ram\`ene \`a $SL_n$ ou
$Sp_{2n-2}$ ; dans le cas de $SO_{2n+1}$ on se ram\`ene \`a $SL_n$ ou
$SO_{2n-1}$ ; dans le cas de $SO_{2n}$ on se ram\`ene \`a $SL_n$ ou
$SO_{2n-2}$ ; dans le cas de $F_4$ on se ram\`ene \`a $Sp_6$ ou
$SO_7$ ; dans le cas de $E_6$ on se ram\`ene \`a $SL_6$ ou
$SO_{10}$ ; dans le cas de $E_7$ on se ram\`ene \`a $SL_7$, $SO_{12}$
ou $E_6$ et dans le cas de $E_8$ on se ram\`ene \`a $SL_8$, $SO_{14}$ ou
$E_7$. 

Ce choix de $P'$ n'est pas toujours possible
pour $SL_{n+1}$, $Sp_{2n}$, $SO_{2n+1}$ et $F_4$ qui n'ont que deux
points extr\^emes. Ceci ne se produit que si $P$ est donn\'e par les
deux points extr\^emes. Dans ce cas, on choisit l'avant dernier point du
diagramme. La vari\'et\'e $R'/R$ est alors $\p^{n-2}\times\pu$ pour les
trois premiers groupes et $\p^2\times\pu$ pour $F_4$. Il existe donc
des courbes lisses sur $R'/R$ d\`es que $n\geq 4$. 

Pour initialiser la r\'ecurrence, il nous faut donc montrer les cas de
$SL_4$, $Sp_4=SO_5$, $Sp_6$ et $SO_7$.

On se ram\`ene, dans le cas de $SL_4$, pour les incidences
point/droite et droite/plan, \`a $SL_3/B$ sur laquelle il y a des
courbes lisses. Il existe une $P'$-orbite de l'incidence point/plan
qui est un fibr\'e vectoriel de degr\'e $(1,1)$ au dessus de
$\pu\times\pu$. Le lemme $4$ (et la remarque $7$) nous permet de
construire des courbes lisses sur cette $P'$-orbite. De m\^eme il
existe une $P'$-orbite de la Grassmannienne des droites qui est un
fibr\'e vectoriel de degr\'e $1$ au dessus de $\p^2$. Le lemme $4$
nous permet de conclure.


Pour $Sp_4$, les vari\'et\'es homog\`enes
associ\'ees \`a des paraboliques maximaux sont $\p^3$ et $Q_3$ la
quadrique de $\p^4$. On va traiter le cas de $Q_3$ en fin de
d\'emonstration et sur $\p^3$ il y a des courbes lisses. Sur $Sp_4/B$
on sait tracer des courbes lisses en se ramenant soit \`a $\p^3$ soit
\`a $Q_3$. 

Pour $Sp_6$ : si le groupe de Picard est de rang $1$, soit le point de
$P$ est le dernier du diagramme alors en prenant le premier point on a
un morphisme vers $Q_3$ qui a des courbes lisses, soit ce n'est pas le
cas et en choisissant le dernier point on arrive dans $\p^2$ avec un
fibr\'e vectoriel de degr\'e $2$ ou $3$. On conclue grace au lemme $3$
(et \`a la remarque $7$). Si le groupe de Picard est de
rang $2$, on a deux cas selon que les points sont aux deux
extr\'emit\'es ou non. Dans le premier cas on arrive, en prenant le
point du milieu, dans $\pu\times\pu$ avec un fibr\'e vectoriel de
degr\'e $(1,1)$ on conclue grace au lemme $3$ (et \`a la remarque
$7$), dans le second on prend un des deux points extr\^emes et on
arrive dans $SL_3/B$ ou $Sp_4/B$.


Pour $SO_7$ : si le groupe
de Picard est de rang $1$, soit le point de $P$ est le dernier du
diagramme alors en prenant le premier point on a un morphisme vers
$\p^3$, soit ce n'est pas le cas et en choisissant le dernier point on
arrive dans $\p^2$ avec un fibr\'e vectoriel de degr\'e $1$ ou $2$. On
conclue grace au lemme $3$ (et \`a la remarque $7$). Si
le groupe de Picard est de rang $2$, on a deux cas selon que les
points sont aux deux extr\'emit\'es ou non. Dans le premier cas on
arrive en prenant le point du milieu dans $\pu\times\pu$ avec un
fibr\'e vectoriel de degr\'e $(1,2)$ on conclue grace au lemme $3$
(et \`a la remarque $7$), dans le second on prend un des
deux points extr\^emes et on arrive dans $SL_3/B$ ou $Sp_4/B$.

LE CAS DE $G_2$ : Le cas de $G_2$ pose un probl\`eme et on ne peut montrer l'existence
de courbes lisses par cette m\'ethode. On va utiliser une seconde
m\'ethode. L'une des vari\'et\'es homog\`enes de $G_2$ est $Q_5$ la
quadrique de $\p^6$ pour laquelle le probl\`eme est d\'ej\`a r\'esolu
il faut donc le v\'erifier pour l'autre vari\'et\'e homog\`ene
associ\'ee \`a un parabolique maximal et l'incidence qui est un
$\pu$-bundle au dessus de chacune des deux.

On a la situation suivante : soient $P_1$ et $P_2$ les paraboliques
maximaux de $G_2$ et $B$ un Borel les contenant. On a alors
l'incidence donn\'ee par les fl\`eches $f:G_2/B \fl G_2/P_1$ et $g :G_2/B
\fl G_2/P_2$. On sait de plus que $f$ est une fibration en $\pu$. On
consid\`ere $G_2/P_2$ plong\'e dans son plongement minimal qui est alors
une quadrique de dimension $5$ (dont on sait qu'elle a des courbes
lisses). On regarde une section hyperplane g\'en\'erale $H$ (sur
laquelle il y a aussi des courbes lisses). On obtient ainsi des
morphismes $f'$ et $g'$ par restriction. Les conditions suivantes sont
r\'ealis\'ees :

\noindent
-Le morphisme $f'$ est birationnel de $g^{-1}((G_2/P_2)\cap H)$ vers
$G_2/P_1$ et c'est un isomorphisme en dehors d'un ferm\'e $E$ de
codimension au moins $1$ de $g^{-1}((G_2/P_2)\cap H)$. Donc $f'(E)$ est
de codimension au moins $2$.

\noindent
- Le morphisme $g$ permet de r\'ealiser $g^{-1}((G_2/P_2)\cap H)-E$ comme
un fibr\'e vectoriel engendr\'e par ses sections au dessus de
$(G_2/P_2)\cap H-Z$ o\`u $Z$ est un ferm\'e de codimension au moins $2$
contenu dans $g(E)$. 

Ainsi, comme il existe des courbes lisses sur $(G_2/P_2)\cap H$ et
que $g'$ est un fibr\'e vectoriel engendr\'e par ses sections, il
existe des courbes lisses sur $g^{-1}((G_2/P_2)\cap H)-E$ et par $f'$
des courbes lisses sur $G_2/P_1$.

LE CAS DE $Q_3$ : Il nous reste \`a montrer qu'il existe des courbes rationnelles lisses
de tous les degr\'es sur $Q_3$. On consid\`ere cette vari\'et\'e comme
les droites isotropes pour une forme symplectique dans un espace de
dimension $4$ (c'est \`a dire comme une vari\'et\'e homog\`ene sous
$Sp_4$). La vari\'et\'e d'incidence avec $\p^3$ est donn\'e par le
fibr\'e projectif associ\'e au faisceau de nulle corellation $E$
d\'efini par la forme symplectique :
$0\fl\oo_{\p^3}(-1)\fl\Omega^1(1)\fl E\fl 0$. On fixe un point $P_0$
de $\p^3$. Soit $H_0$ l'orthogonal de ce point et $L_0$ la droite de
$Q_3$ form\'ee par les droites isotropes de $H_0$. La restriction de
$E$ \`a $H_0$ est donn\'ee par l'extension non triviale
$0\fl\oo_{H_0}\fl E\fl\I_{P_0}\fl 0$. Au dessus de $H_0-P_0$ le
fibr\'e $E$ est une extension non triviale de $\oo_{H_0-P_0}$ par lui
m\^eme qui a une unique section $s$ donn\'ee par $P\mapsto
(P,(PP_0))$. Soit $Z$ l'image de cette section.

On s'int\'eresse maintenant \`a $Q_3-L_0$ et la vari\'et\'e
d'incidence $X$ (qui est au dessus de $\p^3-P_0$). On a une section de
$Q_3-L_0$ vers $X$ donn\'ee par $L\mapsto (L\cap H_0,L)$. Soit $Y$
l'image de cette section, si on reprojette $Y$ vers $\p^3$, son image
est $H_0-P_0$. Plus pr\'ecis\'ement $Y$ est le fibr\'e
$\p_{H_0-P_0}(E)$ priv\'e de la section $Z$. $Y$ est donc un fibr\'e
affine non vectoriel au dessus de $H_0-P_0$, de fibr\'e vectoriel
direction $\oo_{H_0-p_0}$, donn\'e par $\eta\in H^1\oo_{H_0-p_0}$. On
cherche \`a tracer des courbes lisses grace au lemme $4$. Soit $C$ une
courbe rationnelle nodale de degr\'e $d$ sur $H_0-P_0$, on cherche \`a
la relever en une courbe lisse dans $Y$ qui est isomorphe \`a
$Q_3-L_0$. Le lemme $4$ nous dit (qu'il faut et) qu'il suffit que
l'image $\overline{\eta}$ de $\eta$ dans $H^1\oo_C$ soit non nulle (ou
que $C$ soit lisse). Ceci est
vrai d\`es que $d\geq 3$ (si $d\leq 2$ la courbe $C$ est lisse).


EXISTENCE DE COURBES RATIONNELLES PARAM\'ETR\'EES : Pour terminer la
d\'emonstration du th\'eor\`eme $1$, il nous reste \`a prouver la non
vacuit\'e de $\Mor{\a}{G/P}$ pour $\a$ dans le c\^one positif. Si $P$
n'est pas un Borel, en prenant tous les points du diagramme distincts
de ceux de $P$ apr\`es involution, on construit une $P'$-orbite
maximale qui est une tour de fibr\'es affines au dessus de $R'/R$ dont
les fibr\'es vectoriels direction sont engendr\'es par leurs
sections. De plus le choix de $P'$ nous permet de dire que $R'/R$ est
un produit de vari\'et\'es de la forme $G'/B'$ o\`u $B'$ est un Borel
de $G'$. Le lemme $4$ nous permet de nous ramener au cas des Borels.

%

Il reste donc \`a traiter
ce cas. Pour cela on proc\`ede une fois encore par r\'ecurrence sur la
longueur du diagramme de Dynkin. On utilise la construction suivante
en supposant que la longueur du diagramme est au moins $3$. Soit $B$
un Borel de $G$, soit $\a\in A_1(G/B)$ dans le c\^one positif, soient $x$ et $y$ deux points
du diagramme de Dynkin qui ne sont pas sur la m\^eme ar\^ete, et
notons $P_x$ (resp. $P_y$) et $P_{x,y}$ les paraboliques donn\'es par
tous les points du diagramme sauf $x$ (resp. sauf $y$) et par tous les
points du diagramme sauf $x$ et $y$. On a alors le diagramme suivant :
$$\begin{array}{ccc}
 G/B & \stackrel{p'}{\fl} G/P_x\\
 \mapdown{q'} & \mapdown{q}\\
 G/P_y & \stackrel{p}{\fl} G/P_{x,y} \\
\end{array}$$
Toutes les fl\`eches de ce diagramme sont des fibrations en droites
projectives et si $p$ et $q$ sont donn\'es par les fibr\'es $E$ et $F$
de rang $2$, alors $p'$ et $q'$ sont donn\'es par les fibr\'es $q^*E$
et $p^*F$.

Sur $G/P_x$ (resp. $G/P_y$), on sait tracer des courbes de classe
$p'_*\a$ (resp. $q'_*\a$). En effet, la remarque pr\'eliminaire permet
de se ramener au cas des Borels d'un groupe dont le diagramme de
Dynkin est de longueur strictement inf\'erieure \`a celle de $G$ et on
conclue par hypoth\`ese de r\'ecurrence. Ceci nous permet d'affirmer
que dans $G/P_{x,y}$ il existe des courbes v\'erifiant la condition du
lemme $5$ pour le fibr\'e $F$. Ces courbes forment un ouvert $U$ non
vide de $\Mor{p_*p'_*\a}{G/P_{x,y}}$.

De la m\^eme fa\c con on voit que $\Mor{q'_*\a}{G/P_y}$ est non
vide. Son image dans $\Mor{q_*q'_*\a}{G/P_{x,y}}$
($q_*q'_*\a=p_*p'_*\a$) est un ouvert (proposition $4$) non vide. Elle
rencontre donc l'ouvert $U$ pr\'ec\'edent. Ainsi, il existe dans
$\Mor{p_*\a}{G/P_y}$ un \'el\'ement $f:\pu\fl G/P_y$ tel que
$f^*p^*F=\oo_{\pu}\oplus\oo_{\pu}(x)$ avec $x\geq 0$, $x\equiv d\
\![2]$ et $x\leq d$ o\`u $d=p'_*\a\cap T_q$ est le
degr\'e de $p'_*\a$ par rapport \`a la fibration $q$
($T_q$ est le fibr\'e tangent relatif de $q$). Mais le
degr\'e de $\a$ par rapport \`a la fibration $q'$ est $\a\cap
T_{q'}=\a\cap {p'}^*T_q={p'}^*(p'_*\a\cap T_q)=d$. Ainsi $f$ v\'erifie les
conditions du lemme $5$ pour $p^*F$ et ceci nous permet de construire
un rel\`evement dans $G/B$.

On s'est ainsi ramen\'e aux cas des groupes dont le diagramme est de
longueur au plus $2$. Pour $SL_2$ et $SO_4=SL2\times SL_2$ c'est
\'evident. Il nous reste donc \`a montrer le cas vari\'et\'es de
drapeaux complets des groupes $SL_3$, $SO_5=Sp_4$ et $G_2$. On
proc\`ede pour les trois de la m\^eme fa\c con et on note $G$ l'un de
ces trois groupes. Soit $X$ la vari\'et\'es homog\`ene qui correspond
aux droites (pour le second groupe on le consid\`ere comme $Sp_4$) et
$Y$ l'autre vari\'et\'e homog\`ene (qui correspond aux points). Soit
$B$ un Borel de $G$, soit $\vp$ la fibration en droites
projectives $G/B\fl X$ et de fibr\'e tangent relatif $T_{\vp}$ et soit
$\a\in A_1(G/B)$ de degr\'es $d_1$ et $d_2$ par rapport \`a $Y$ et $X$,
alors $\a\cap T_{\vp}=2d_1-d_2$. Soit $f:\pu\fl X$ un morphisme de
degr\'e $d_2$. Soit $E$ la restriction du fibr\'e tautologique de la
Grassmannienne $\G(2,m)$ ($m=3,4,7$ selon les cas) \`a $X$, c'est un
fibr\'e vectoriel associ\'e \`a la fibration $\vp$. On a
n\'ecessairement $f^*E=\oo_{\pu}(a)\oplus\oo_{\pu}(b)$ avec $0\leq
a\leq b$ et $a+b=d_2$. Ainsi $f$ v\'erifie les hypoth\`eses du lemme
$5$ d\`es que $2d_1-d_2\geq d_2\geq b-a$ ie $d_1\geq d_2$ ce qui dans ce cas
nous permet de relever $f$ dans $G/B$. 

Il nous reste donc \`a tracer des courbes sur $G/B$ pour $d_2\geq
d_1$. Pour ceci on trace une courbe de bidegr\'e $(d_1,d'_2)$ dans
$G/B$ avec $d_1\geq d'_2$. Ceci nous permet de dire qu'il existe un
faisceau $F$ associ\'e \`a la fibration de $G/B$ au dessus de $Y$ tel
que si $f:\pu\fl Y$ est le morphisme induit, alors $f^*F$ v\'erifie
les hypoth\`ese du lemme $5$ (le degr\'e relatif $d$ est alors
$2d'_2-d_1$ pour $SL_3$, $2d'_2-2d_1$ pour $Sp_4$ et $2d'_2-3d_1$ pour
$G_2$). On sait ainsi que $f^*F=\oo_{\pu}\oplus\oo_{\pu}(x)$
avec $x\geq 0$, $d\equiv x\ \![2]$ et $d\geq x$. Mais alors si on veut une courbe de
bidegr\'e $(d_1,d_2)$ avec $d_2\geq d'_2$ quelconque, il suffit de
relever $f$ dans $G/B$ ce qui est possible car le degr\'e relatif est
alors plus grand que $d$ (et donc plus grand que $x$) et que la
parit\'e ne change pas.




\vs 0.4 cm

\noi
{\bf Remarque 8} : Soit $C$ un courbe dont la classe dans $A_1(G/P)$ est dans le
c\^one positif mais pas dans le c\^one strictement positif. On
distingue deux types de points parmis les points d\'efinissant $P$
dans le diagramme de Dynkin, ceux pour lesquels le degr\'e de $C$ est
strictement positif et ceux pour lesquels le degr\'e de $C$ est
nul. On note $P'$ le parabolique obtenu en prenant les points du
deuxi\`eme type. On voit que $P\subset P'$ et on a ainsi un morphisme
$G/P\fl\gsp'$. La courbe $C$ est trac\'ee dans une fibre de ce
morphisme. Si on regarde les sommets du diagramme de Dynkin de $P$
dans les composantes connexes du diagramme de Dynkin de $G$ priv\'e de
$P'$, le produit de vari\'et\'es homog\`enes ainsi d\'efini est
isomorphe \`a la fibre. Pour savoir si il existe des courbes lisses,
on est ainsi ramen\'e aux c\^ones strictements positifs de ce produit
de vari\'et\'es homog\`enes.

\vs 0.5 cm

\centerline{\bf\uppercase{{Une d\'esingularisation des vari\'et\'es de
Schubert}}}

\vs 0.5 cm

\setcounter{section}{0}

On va proposer dans cette partie une application des $P'$-orbites
d\'efinies dans la partie pr\'ec\'edente. En s'inspirant en grande partie
de la d\'esingularisation des vari\'et\'es de Schubert donn\'ee par
M. Demazure [D], on construit une d\'esingularisation plus fine
des vari\'et\'es de Schubert. 

Soit $B$ un Borel et $w\in W$, dans son article [D], M. Demazure
construit une suite de paraboliques $Q_i$ \textit{minimaux}
(associ\'es \`a un unique sommet du diagramme) tels que
$wQ_i\subset\overline{BwB}$ (ce qui revient \`a dire
$\gq_i\subset\gb+\gb'$), $Q_1$ contient $B$, $Q_i\cap Q_{i+1}$
contient un Borel, la suite des $Q_i\cap B$ est d\'ecroissante et
$\sum\gq_i=\gb+\gb'$ (on a not\'e $B'=w^{-1}Bw$). La
d\'esingularisation est alors donn\'ee par le quotient par $B$ de la
vari\'et\'e $X'=Y'/G''$ o\`u $Y'=\prod_iQ_i$ et $G''=\prod_i(Q_i\cap
Q_{i+1})$. Le produit dans $G$ donne un morphisme de $Y'$ dans
$\overline{BwB}$ qui est invariant sous $G''$ et par passage au
quotient sous $B$ on a le morphisme $D:X'\fl\overline{BwB/B}$.

Lorsque $Y'$ est le produit d'un ou deux groupes, le morphisme de
$X'$ dans $G$ est bijectif sur son image. Lorsque $Y'$ a plus de
facteurs, le morphisme de $X'$ dans $\overline{BwB}$ est seulement
birationnel. Plus $Y'$ a de facteurs, plus le morphisme est
succeptible d'avoir des contractions. Ce que l'on va faire ici est
r\'eduire, de fa\c con canonique, le nombre de facteurs de $Y'$. On va
ainsi regrouper les paraboliques $Q_i$ pour les remplacer par des
paraboliques plus grand. On cherche donc des paraboliques $P_i$
\textit{les plus grand possible} tels que $wP_i\subset\overline{BwB}$,
$P_1$ contient $B$, $P_i\cap P_{i+1}$ contient un Borel, la suite des
$P_i\cap B$ est d\'ecroissante et $\sum\gp_i=\gb+\gb'$.

Cette consruction nous permettra de donner un crit\`ere pour qu'une
vari\'et\'e de Schubert soit une $P''$-orbite pour un parabolique
$P''$ de $G$. On donnera ainsi une condition suffisante
(non n\'ecessaire) de lissit\'e des vari\'et\'e de Schubert.

On utilise l'abus de notation suivant : si $\a$ est une racine et
$\gp$ l'alg\`ebre de Lie d'un parabolique, on dit que $\a\in\gp$ si
$\a$ est une valeur propre pour l'action du tore $\gh$ sur $\gp$ ou
encore si $\gg_{\a}$ apparait dans la d\'ecomposition de $\gp$. En
d'autres termes, on identifie l'alg\`ebre de Lie $\gp$ et la partie
parabolique de $P$ qui est l'ensemble des racines de $P$.

\section{Construction des paraboliques}

On commence par montrer (lemme $1$) qu'il existe un parabolique $P_1$
contenu dans $w^{-1}\overline{BwB}$, contenant $B$ et qui est maximal
pour cette propri\'et\'e. Ceci nous permet de
construire par r\'ecurrence une suite de paraboliques plus gros que
ceux de M. Demazure qui donnent la d\'esingularisation annonc\'ee.




%
%

\vs 0.4 cm

\noi
{\bf Lemme 1} : \textit{Soient $\gb$ et $\gb'$ deux Borels de $\gg$
contenant le tore $\gh$. Il existe un unique parabolique $\gp_1$ et un
unique parabolique $\gp'_1$ qui sont maximaux pour les propri\'et\'es
suivantes : $\gb\subset\gp_1\subset\gb +\gb'$ et
$\gb'\subset\gp'_1\subset\gb +\gb'$.}

\vs 0.2 cm

\dm :
Il suffit de montrer l'existence et l'unicit\'e de $\gp_1$, par
sym\'etrie celle de $\gp'_1$ en d\'ecoulera. On prend alors pour
$\gp_1$ le parabolique contenant $\gb$ d\'efini par les oppos\'ees des
racines simples de $\gb$ qui sont dant $-\gb'$. V\'erifions que ce
parabolique est contenu dans $\gb+\gb'$. Soit $\a\in\gp$. Si
$\a\in\gb$, on a termin\'e. Sinon, $\a$ s'\'ecrit $\sum(-\a_i)$ o\`u
$\a_i$ est une racine simple de $\gb$ telle que $-\a_i\in\gb'$. Mais
alors $\sum(-\a_i)\in\gb'$ donc $\a\in\gb'$.

Soit maintenant $\gp_2$ un parabolique v\'erifiant ces
conditions. Soit $\a$ une racine simple de $\gb$ telle que
$\-a\in\gp_2$. Alors, $-\a\in\gb'$ et donc $-\a\in\gp_1$. Ainsi
$\gp_2\subset\gp_1$.

\vs 0.4 cm

Ce lemme permet de construire une suite de paraboliques qui serviront
\`a la d\'esingularisation. Soient $B$ et $B'$ deux Borels.
On note $P_1$ et $P'_1$ les paraboliques construits \`a partir de $B$
et $B'$ et du lemme $1$.
On construit ainsi par r\'ecurrence deux suites de Borels $B_n$ et
$B'_n$ ($B_1=B$ et $B'_1=B'$) et deux suites de paraboliques $P_n$ et
$P'_n$ tels que $B_n$ est contenu dans $P_{n-1}$ et $P_n$ et de m\^eme
$B'_n$ est contenu dans $P'_{n-1}$ et $P'_n$. En effet, supposons $B_n$,
$B'_n$, $P_n$ et $P'_n$ construits, alors on construit $B_{n+1}$ (et
par sym\'etrie $B'_{n+1}$) de la fa\c con suivante : on d\'ecrit les
$\a\in\gp_n$ qui sont dans $\gb_{n+1}$ : si $\a\in\gp_n$ est telle que
$-\a\not\in\gp_n$ alors $\a\in\gb_{n+1}$. Si $\a\in\gp_n\cap\gp'_n$ et
$-\a\in\gp_n$ mais $-\a\not\in\gp'_n$ alors $\a\in\gb_{n+1}$. Si
$\a\in\gp_n$ mais $\a\not\in\gp'_n$ et $-\a\in\gp_n\cap\gp'_n$ alors
$\a\not\in\gb_{n+1}$. Enfin, si $\a\in\gp_n\cap\gp'_n$ et
$-\a\in\gp_n\cap\gp'_n$ alors $\a\in\gb_{n+1}\Leftrightarrow\a\in\gb_n$. Une
fois les Borels $B_{n+1}$ et $B'_{n+1}$ d\'efinis, on d\'efinit
$P_{n+1}$ et $P'_{n+1}$ comme \'etant les paraboliques obtenus \`a
partir de $B_{n+1}$ et $B'_{n+1}$ et du lemme $1$.

Pour le lemme suivant, on a besoin du :

\vs 0.4 cm

\noi
{\bf Fait 1} : \textit{Soit $\gp$ un parabolique et soient $\a\in\gp$
et $\a'\not\in\gp$ deux racines, on a l'implication :}
$$\a+\a'\in\gp\Rightarrow-\a\not\in\gp$$


\dm :
On va en fait montrer que la premi\`ere condition implique la
condition suivante (qui est \'equivalente \`a celle que l'on
cherche) : $\a$ appartient \`a tous les Borels de
$\gp$.

Soit $\gb$ un
Borel de $\gp$, $\a'$ s'\'ecrit $\sum(-\a'_j)$ o\`u $\a'_j$ est une
racine simple de $\gb$. Comme $\a'\not\in\gp$, il existe au moins un
$j$ pour lequel $-\a'_j\not\in\gp$. Supposons que $\a\not\in\gb$, alors $\a$
s'\'ecrit $\sum(-\a_i)$ o\`u $\a_i$ est une racine simple de
$\gb$. Comme $\a\in\gp$, on a pour tout $i$ : $-\a_i\in\gp$. Mais
alors $\a+\a'\in\gp$ impose que $\a+\a'$ s'\'ecrive $\sum(-\a''_k)$,
o\`u pour tout $k$ la racine $\a''_k$ est simple dans $\gb$ et telle
que $-\a''_k\in\gp$. Mais alors les $\a'_j$ sont contenus dans les
$\a''_k$ ce qui impose que pour tout $j$ on ait $-\a'_j\in\gp$ ce qui
est une contradiction. 

\vs 0.4 cm

\noi
{\bf Lemme 2} : \textit{Les parties $\gb_{n+1}$ et $\gb'_{n+1}$
sont les alg\`ebres de Lie de Borels et on a pour tout $n$ :
$\gb_{n+1}+\gb'_{n+1}\subset\gb_n+\gb'_n$ et
$\gb_n\cap\gb'_1\subset\gb_{n+1}\cap\gb'_1\subset\gb'_{n+1}\cap\gb'_1
\subset\gb'_n\cap\gb'_1$.}

\vs 0.2 cm

\dm :
On commence par montrer que les parties $\gb_{n+1}$ et $\gb'_{n+1}$
sont les alg\`ebres de Lie de Borels. Il suffit par sym\'etrie de le
faire pour $\gb_{n+1}$. On doit donc montrer que $\gb_{n+1}$ est
stable, que pour tout $\a\in\gg$, l'une des deux racines $\a$ ou $-\a$
est dans $\gb_{n+1}$ et que l'on a jamais les deux en m\^eme
temps. On commence par la stabilit\'e. Soit $\a$ et $\a'$ des racines
de $\gb_{n+1}$. 

\noi
$\bullet$ Si $-\a$ et $-\a'$ ne sont pas dans $\gp_{n+1}$, alors
c'est aussi le cas de $-(\a+\a')$ et donc $\a+\a'\in\gb_{n+1}$. 

\noi
$\bullet$ Si $-\a$
n'est pas dans $\gp_n$ et que l'on a $\a'\in\gp_n\cap\gp'_n$ et
$-\a'\in\gp_n$ mais $-\a'\not\in\gp'_n$, alors si $-(\a+\a')\in\gp_n$,
on a (Fait $1$) $\a'\not\in\gp_n$ ce qui est absurde. Donc
$-(\a+\a')\not\in\gp_n$ et donc $\a+\a'\in\gb_{n+1}$. 

\noi
$\bullet$ Si $-\a$
n'est pas dans $\gp_n$ et que l'on a $\a'\in\gp_n\cap\gp'_n$ et
$-\a'\in\gp_n\cap\gp'_n$, alors de la m\^eme fa\c con on doit avoir
$-(\a+\a')\not\in\gp_n$ et donc $\a+\a'\in\gb_{n+1}$. 

\noi
$\bullet$ Si
$\{\a,\a'\}\subset\gp_n\cap\gp'_n$ et $\{-\a,-\a'\}\subset\gp_n$ mais
$-\a\not\in\gp'_n$, $-\a'\not\in\gp'$ alors
$\a+\a'\in\gp_n\cap\gp'_n$, $-(\a+\a')\in\gp_n$ et
$-(\a+\a')\not\in\gp'_n$ donc $\a+\a'\in\gb_{n+1}$. 

\noi
$\bullet$ Si
$\a\in\gp_n\cap\gp'_n$ et $-\a\in\gp_n$ mais $-\a\not\in\gp'_n$ et que
l'on a $\a'\in\gp_n\cap\gp'_n$ et $-\a'\in\gp_n\cap\gp'_n$, alors
$\a+\a'\in\gp_n\cap\gp'_n$ et $-(\a+\a')\in\gp_n$. Si de plus
$-(\a+\a')\in\gp'_n$, on a (Fait $1$) $\a'\not\in\gp'_n$ ce qui est
absurde donc $-(\a+\a')\not\in\gp'_n$ et donc $\a+\a'\in\gb_{n+1}$. 

\noi
$\bullet$ Enfin, si $\{\a,\a'\}\subset\gp_n\cap\gp'_n$,
$\{-\a,-\a'\}\subset\gp_n\cap\gp'_n$ et $\{\a,\a'\}\subset\gb_n$, alors
$\a+\a'\in\gp_n\cap\gp'_n$, $-(\a+\a')\in\gp_n\cap\gp'_n$ et
$\a+\a'\in\gb_n$ donc $\a+\a'\in\gb_{n+1}$.

Soit maintenant $\a\in\gg$. Si $\a\not\in\gp_n$ alors $-\a\in\gp_n$
et donc $-\a\in\gb_{n+1}$. De m\^eme si $-\a\not\in\gp_n$ alors $\a\in\gp_n$ et donc
$\a\in\gb_{n+1}$. Il reste donc les racines $\a\in\gp_n$ telles que
$-\a\in\gp_n$. On sait que $\a$ ou $-\a$ est dans $\gp'_n$, on peut
donc supposer (quitter \`a \'echanger $\a$ et $-\a$) que
$\a\in\gp'_n$. On a alors deux cas : $-\a\not\in\gp'_n$ ou
$\a\in\gp'_n$. Dans le premier cas on sait que $\a\in\gb_{n+1}$, dans
le second on a $\a\in\gb_{n+1}\Leftrightarrow\a\in\gb_n$. Or on sait que $\gb_n$
est un Borel donc $\a$ ou $-\a$ est dans $\gb_n$ et ainsi $\a$ ou
$-\a$ est dans $\gb_{n+1}$. Il reste \`a voir que l'on a pas $\a$ et
$-\a$ dans $\gb_{n+1}$. Si c'est le cas on sait que $\a$ et $-\a$ sont
dans $\gp_n$. Si $\a\not\in\gp'_n$ alors $-\a\in\gp'_n$ et ceci impose
que $\a\not\in\gb_{n+1}$ ce qui est absurde. Par sym\'etrie on peut
donc supposer que $\a$ et $-\a$ sont dans $\gp'_n$, mais alors
$\a\in\gb_{n+1}\Leftrightarrow\a\in\gb_n$ et comme $\gb_n$ est un Borel on ne
peut avoir $\a$ et $-\a$ dans $\gb_{n+1}$. 

Par construction on sait que $\gb_{n+1}\subset\gp_n$ et
$\gb'_{n+1}\subset\gp_n$ ce qui nous donne que
$\gb_{n+1}+\gb'_{n+1}\subset\gp_n+\gp'_n=\gb_n+\gb'_n$.

On proc\`ede par r\'ecurrence en supposant que
$\gb_n\cap\gb'_1\subset\gb'_n\cap\gb'_1$ cette propri\'et\'e \'etant
\'evidement vraie pour $n=1$. 
On a $\gb'_{n+1}\cap\gb'_1\subset
(\gb_n+\gb'_n)\cap\gb'_1\subset\gb'_n\cap\gb'_1$ par hypoth\`ese de
r\'ecurrence. 

De m\^eme, on a $\gb_{n+1}\cap\gb'_1\subset\gb'_n\cap\gb'_1$. Soit
alors $\a\in\gb_{n+1}\cap\gb'_1$. On sait alors que
$\a\in\gp_n\cap\gb'_n\subset\gp_n\cap\gp'_n$ et on a les cas suivants
:

\noi
$\bullet$ $-\a\not\in\gp'_n$ alors $\a$ est dans tous les
Borels de $\gp'_n$ et donc $\a\in\gb'_{n+1}$.

\noindent
$\bullet$ $-\a\in\gp'_n$ mais $-\a\not\in\gp_n$ alors $-\a\in\gb'_n$ (sinon
$-\a\not\in\gb_n+\gb'_n$ alors que $-\a\in\gp'_n$) et donc
$\a\not\in\gb'_n$ et ce cas ne se produit pas.

\noindent
$\bullet$ $-\a\in\gp_n\cap\gp'_n$ alors on a $\a\in\gb'_n$ et donc
$\a\in\gb'_{n+1}$.

\noindent
On conclue ainsi que $\gb_{n+1}\cap
\gb'_1\subset\gb'_{n+1}\cap\gb'_1$. 

\vs 0.1 cm

On sait que $\gb_n\cap\gb'_1\subset\gb'_n\cap\gb'_1$. Soit
$\a\in\gb_n\cap\gb'_1$. Si $-\a\not\in\gp_n$ alors $\a$ est dans tous
les Borels de $\gp_n$ et donc $\a\in\gb_{n+1}$. Supposons
$-\a\in\gp_n$. Comme $\a\in\gb'_n$ alors $\a\in\gp'_n$. Mais alors on
a les deux cas suivants :

\noindent
$\bullet$ $-\a\not\in\gp'_n$ alors $-\a\in\gb_n$ (sinon
$-\a\not\in\gb_n+\gb'_n$ alors que $-\a\in\gp_n$) ce qui impose
$\a\not\in\gb_n$ et ce cas est exclu.

\noindent
$\bullet$ $-\a\in\gp'_n$ et on a $\a\in\gp_n\cap\gp'_n$ et
$-\a\in\gp_n\cap\gp'_n$ donc comme $\a\in\gb_n$ on $\a\in\gb_{n+1}$.

\noindent
On conclue ainsi que $\gb_{n}\cap
\gb'_1\subset\gb_{n+1}\cap\gb'_1$.

\vs 0.4 cm

On a ainsi construit deux suites de Borels et deux suites de
paraboliques. On s'arr\^ete d\`es que $P_n\cap P'_n$ contient un
Borel (il est n\'ecessaire d'avoir cette propri\'et\'e pour que le
morphisme $\pi$ que l'on construit dans la suite et qui est la
d\'esingularisation soit propre). On a le lemme suivant qui nous
permet de dire que notre construction s'arr\^ete.

\vs 0.4 cm

\noi
{\bf Lemme 3} : \textit{Si $\gp_n\cap\gp'_n$ ne contient pas de Borel
alors l'inclusion $\gb_{n+1}+\gb'_{n+1}\subset\gb_n+\gb'_n$ est
stricte.}

\vs 0.2 cm

\dm :
Supposons que $\gp_n\cap\gp'_n$ ne contient pas de Borel et que
$\gb_{n+1}+\gb'_{n+1}=\gb_n+\gb'_n$. Si il existe $\a\in\gp_n$ telle
que $\a\not\in\gp'_n$ et que $-\a\in\gp_n\cap\gp'_n$ alors
$\a\in\gb_n+\gb'_n$ et $-\a\in\gb_{n+1}\cap\gb'_{n+1}$ donc
$\a\not\in\gb_{n+1}+\gb'_{n+1}$ ce qui est impossible. De m\^eme si il
existe $\a\in\gp_n\cap\gp'_n$ telle que $-\a\in\gp'_n$ mais
$-\a\not\in\gp_n$ alors $-\a\in\gb_n+\gb'_n$ et
$-\a\not\in\gb_{n+1}+\gb'_{n+1}$ ce qui est impossible. 
Les racines de $\gb_n$ sont donc d'un des trois types suivants :
$\a\in\gp_n\cap\gp'_n$ et $-\a\in\gp_n\cap\gp'_n$ ou $\a\in\gp_n$,
$\a\not\in\gp'_n$, $-\a\in\gp'_n$ et $-\a\not\in\gp_n$ ou
$\a\in\gp_n\cap\gp'_n$ et $-\a\not\in\gp_n+\gp'_n$.

Soit maintenant $\gb$ un Borel de $\gp_n$ tel que
${\rm{Card}}(\gb\cap\gp_n\cap\gp'_n)$ est maximal (ici on appelle
${\rm{Card}}(\gp)$ le nombre de racines qui apparaissent dans
$\gp$). Alors il existe une racine simple $\a$ de $\gb$ telle que
$\a\not\in\gp_n\cap\gp'_n$ et $-\a\not\in\gp_n\cap\gp'_n$ (pour $\a$
c'est clair sinon $\gb$ serait contenu dans $\gp_n\cap\gp'_n$, si
$-\a\in\gp_n\cap\gp'_n$, alors $s_{\a}(\gb)$ est un Borel de $\gp_n$
tel que
$s_{\a}(\gb)\cap\gp_n\cap\gp'_n=(\gb\cap\gp_n\cap\gp'_n)\cup\{-\a\}$
ce qui contredit la maximalit\'e). On sait alors que $\a\in\gp_n$ mais
$\a\not\in\gp'_n$ et $-\a\in\gp'_n$ mais $-\a\not\in\gp_n$. On montre
maintenant que $-\a\in\gp_n$ ce qui sera une contradiction. Pour cela
il suffit de montrer que pour tout $\b\in\gb_n$ on a
$\b-\a\in\gb_n+\gb'_n$. Mais la remarque faite au debut nous permet de
dire que l'on les trois cas suivants :

\noi
$\bullet$
$\b\in\gp_n\cap\gp'_n$ et $-\b\in\gp_n\cap\gp'_n$, alors $\b\in\gp'_n$
et $-\a\in\gp'_n$ donc $\b-\a\in\gp'_n\subset\gb_n+\gb'_n$.

\noi
$\bullet$
$\b\in\gp_n\cap\gp'_n$ et $-\b\not\in\gp_n+\gp'_n$, alors $\b\in\gp'_n$
et $-\a\in\gp'_n$ donc $\b-\a\in\gp'_n\subset\gb_n+\gb'_n$.

\noi
$\bullet$
$\b\in\gp_n$ mais $\b\not\in\gp'_n$ et $-\b\in\gp'_n$ mais
$-\b\not\in\gp_n$, alors $\b$ est dans tous les Borels de $\gp_n$ et
en particulier dans $\gb$. On peut donc \'ecrire $\b=\sum\a_i$ o\`u
les $\a_i$ sont des racines simples de $\gb$. Mais alors
$\b-\a=\sum\a_i-\a$ est une racine de $\gg$ si et seulement si
$\a\in\{\a_i\}$ (car $\a$ est une racine simple) et donc
$\b-\a\in\gb\subset\gp_n\subset\gb_n+\gb'_n$.

\section{La d\'esingularisation}

Soient $P$ et $P'$ deux paraboliques contenant un m\^eme Borel et soit
$w\in W$. On peut maintenant construire la d\'esingularisation de la
vari\'et\'e $\overline{P'wP/P}$ qui est l'adh\'erence de la
$P'$-orbite $P'wP/P$ dans $\gsp$. Le lemme $1$ de la partie
pr\'ec\'edente nous permet de construire un Borel $B\subset P$ et
$w'\in\overline{w}$ tels que $\gb+\gb'=\gp+\gp'$ ($\gb'=w'(\gb)$). On
peut donc se contenter de construire la suite de paraboliques pour $B$
et $w'$. La d\'esingularisation de $\overline{P'wP/P}$ sera alors le
quotient de celle de $\overline{Bw'B/B}$ par $P$. Ce sera une fibration
en $P/B$. On note
$B'=w'(B)$. Pour cette construction, on s'inspire directement de celle
de Demazure [D]. Les lemmes pr\'ec\'edents nous ont permis de
construire des suites de paraboliques et de Borels. De plus le lemme
$3$ nous permet de dire qu'a partir d'un certain rang
$\gp_n\cap\gp'_n$ contiendra un Borel : tant que ce n'est pas le cas
la suite des $\gb_n+\gb'_n$ est strictement d\'ecroissante (en
dimension) et sera donc constante \`a partir d'un certain rang.

On note $P_i$, $P'_i$, $B_i$ et $B'_i$ les paraboliques associ\'es aux
alg\`ebres de Lies $\gp_i$, $\gp'_i$, $\gb_i$ et $\gb'_i$. On construit
la vari\'et\'e suivante : $X=P'_1\times^{P'_1\cap P'_2}P'_2\cdots
P_2\times^{P_1\cap P_2}P_1$ qui est le quotient de $Y=P'_1\times
P'_2\cdots P_2\times P_1$ par $G'=(P'_1\cap P'_2)\times \cdots \times
(P_1\cap P_2)$. On a un morphisme de $Y$ vers $G$ donn\'e par le
produit (et par multiplication par $w'$ pour arriver dans
$\overline{Bw'B}$) qui est invariant sous l'action de $G'$ ce qui nous donne un
morphisme de $X$ dans $G$. On voit alors que $B$ (et m\^eme $P$)
agit \`a droite sur ces vari\'et\'es et on obtient ainsi un morphisme
$\pi$ de $X/B$ (ou $X/P$) vers $G/B$ (ou $G/P$).

On va comparer cette construction avec celle donn\'ee par M. Demazure
dans [D] : on a, pour former $P_1$, regroup\'e les paraboliques que
choisit M. Demazure. On voit ainsi que $\pi$ factorise la
d\'esingularisation de Demazure (cf. proposition $1$).

L'image de $\pi$ est l'adh\'erence de la cellule $BwB/B$ et $X/B$ est
le quotient de $P'_1\times P'_2\times\cdots\times P_2\times P_1/B$ par
$G'$ et est donc lisse. On peut aussi le voir en consid\'erant la
filtration suivante :
$$X/B\fl P'_1\times^{P'_1\cap P'_2}P'_2\times\cdots
\times P_2/(P_1\cap P_2)\fl\cdots\fl
P'_1\times^{P'_1\cap P'_2}P'_2/(P'_2\cap P'_3)\fl P'_1/(P'_1\cap P'_2)$$
o\`u la premi\`ere fl\`eche est une fibration en $P_1/B$,la seconde
en $P_2/(P_1\cap P_2)$ et ainsi de suite chacune des fl\`eches est
un fibration en $P_{n+1}/(P_n\cap P_{n+1})$ ou en $P'_n/(P'_n\cap
P'_{n+1})$. Ce qui prouve que $X/P$ est lisse.

\vs 0.4 cm

\noi
{\bf Remarque 1} : On va voir que cette d\'esingularisation est plus
fine que celle de [D] : la d\'esingularisation de Demazure se
factorise par $\pi$. Notre d\'esingularisation est alors un
isomorphisme sur un ouvert plus grand. En fait notre
d\'esingularisation est bijective sur tout l'ouvert $w'P'_1B/B$ qui
contient $w'B'B/B=Bw'B/B$.

\vs 0.4 cm

\noi
{\bf Exemple 1} : Si $Bw'B/B$ est la cellule maximale de $G/B$
alors ceci signifie que $\gb+\gb'=\gg$ et donc
$\gp_1=\gp'_1=\gg$. Ainsi notre d\'esingularisation est $G/B$ (qui
\'etait d\'ej\`a lisse) alors que celle de [D] \'etait plus
compliqu\'ee et notamment pas un isomorphisme.

\vs 0.4 cm

Pour montrer la factorisation annonc\'ee, on construits des
paraboliques \textit{minimaux} contenus dans les $P_i$ et qui redonne
la d\'esingularisation de M. Demazure.

\vs 0.4 cm

\noi
{\bf Lemme 4} : \textit{Soit $n$ le plus petit entier tel que
$\gp_n\cap\gp'_n$ contienne un Borel. Il existe $\gb_{n+1}$ un Borel
de $\gp_n\cap\gp'_n$ tel que
$\gb_n\cap\gb'_1\subset\gb_{n+1}\cap\gb'_1\subset\gb'_n\cap\gb'_1$.}


\vs 0.2 cm

\dm :
On modifie $\gb_n$ pour le faire appartenir \`a $\gp'_n$ : si
$\gb_n\subset\gp'_n$, on pose $\gb_{n+1}=\gb_n$. Sinon, il existe $\a$
racine simple de $\gb_n$ telle que $\a\not\in\gp'_n$ mais
$-\a\in\gp_n$. En effet, soit $\a$ une racine simple de $\gb_n$ qui
n'est pas dans $\gp'_n$ (qui existe car $\gb_n\not\subset\gp'_n$),
alors si $-\a\not\in\gp_n$, on a $\a\not\in\gp_n\cap\gp'_n$ et
$-\a\not\in\gp_n\cap\gp'_n$ ce qui contredit le fait que
$\gp_n\cap\gp'_n$ contient un Borel. On regarde alors
$s_{\a}(\gb_n)$. On a $s_{\a}(\gb_n)\subset\gp_n$ et
$s_{\a}(\gb_n)\cap\gp'_n=(\gb_n\cap\gp'_n)\cup\{-\a\}$. De plus, si
$\a\in\gb'_1$, alors
$\a\in\gb_n\cap\gp'_1\subset\gb'_n\cap\gb'_1\subset\gp'_n$ ce qui est
absurde. Donc $\a\not\in\gb'_1$ et
$s_{\a}(\gb_n)\cap\gb'_1=(\gb_n\cap\gb'_1)\cup\{-\a\}$. On continue ce
processus tant que le Borel n'est pas contenu dans $\gp'_n$. On
obtient ainsi le Borel $\gb_{n+1}$ qui v\'erifie
$\gb_{n+1}\subset\gp_n\cap\gp'_n$,
$\gb_n\cap\gb'_1\subset\gb_{n+1}\cap\gb'_1$. Enfin,
$\gb_{n+1}\cap\gb'_1\subset(\gb_n+\gb'_n)\cap\gb'_1\subset\gb'_n\cap\gb'_1$.

\vs 0.4 cm

\noi
{\bf Proposition 1} : \textit{Soit $S$ une vari\'et\'e de Schubert, il
existe une d\'esingularisation de Demazure $D:X'\fl S$ et un morphisme
$X'\fl X$ qui s'ins\`ere dans le diagramme suivant :}
$$\begin{array}{ccc}
 X' & \stackrel{D}{\fl} & S\\
 \mapdown{} & & \mapdown{{\rm{Id}}}\\
 X & \stackrel{\pi}{\fl} & S\\
\end{array}$$

\vs 0.2 cm 

\dm :
On commence par le lemme suivant :

\vs 0.1 cm

\noi
{\bf Lemme 5} : \textit{Soit $\gp$ un parabolique, $\gb$ un Borel,
$\gb'$ et $\gb''$ des Borels de $\gp$ tels que
$\gb\cap\gb'\subset\gb\cap\gb''$, alors il existe une suite de Borels
$(\gb_i)_{1\leq i\leq n}$ de $\gp$ tels que $\gb_1=\gb'$,
$\gb_n=\gb''$, $\gb_i\cap\gb\subset\gb_{i+1}\cap\gb$ et ${\rm{Card}}(\gb_{i+1}\cap\gb)={\rm{Card}}(\gb_i\cap\gb)+1$.}

\vs 0.2 cm

\dm :
On peut se placer dans le cas o\`u l'inclusion
$\gb\cap\gb'\subset\gb\cap\gb''$ est stricte. Alors il existe $\a$ une
racine simple de $\gb'$ telle que $-\a\in\gb\cap\gb''$. En effet,
sinon toute les racines simples de $\gb'$ sont telles que
$-\a\not\in\gb\cap\gb''$ c'est \`a dire $\a\in\gb$ ou
$\a\in\gb''$. Mais si $\a\in\gb$ alors
$\a\in\gb\cap\gb'\subset\gb\cap\gb''$ et donc dans tous les cas
$\a\in\gb''$. Ceci impose que $\gb'\subset\gb''$ et donc $\gb'=\gb''$
ce qui est impossible en raison de l'inclusion stricte.

On pose alors $\gb_2=s_{\a}(\gb')$ et on a ($-\a\in\gb''$ donc
$-\a\in\gp$) $\gb_2\subset\gp$ et
$\gb_2\cap\gb=(\gb'\cap\gb)\cup\{-\a\}\subset\gb\cap\gb''$. On
recommence le processus tant que l'inclusion de
$\gb_i\cap\gb\subset\gb\cap\gb''$ est stricte.

\vs 0.2 cm

Ce lemme nous permet de construire pour tout $1\leq k\leq n$ deux suites
de Borels $(\gb_{k,i})_{1\leq i\leq r_k}$ (resp. $(\gb'_{k,i})_{1\leq
i\leq r'_k}$) de $\gp_i$ (resp. $\gp'_i$) tels que les
$\gb_{k,i}\cap\gb'_1$ (resp. $\gb'_{k,i}\cap\gb'_1$) forment une suite
croissante (resp. d\'ecroissante) pour l'ordre lexicographique et que
leur dimension augente (resp. descende) de exactement $1$ \`a
chaque. De m\^eme le lemme $5$ nous permet de construire deux suites
de Borels $(\gb_{n+1,i})_{1\leq i\leq r_{n+1}}$ de $\gp_n$ (entre
$\gb_n$ et $\gb_{n+1}$) et $(\gb'_{n+1,i})_{1\leq i\leq r'_{n+1}}$ de
$\gp'_n$ (entre $\gb_{n+1}$ et $\gb'_{n}$)
telles que les dimensions de leurs intersections avec $\gb'_1$
forment une suite croissante pour la premi\`ere et d\'ecroissante
pour la seconde. Ces deux suites compl\`etent les deux premi\`eres
suites en une seule qui est telle que les intersections avec $\gb'_1$
forment une suite dont les dimensions prennent une seule fois toutes
les valeurs entre ${\rm{dim}}(\gb'_1)$ et
${\rm{dim}}(\gb_1\cap\gb'_1)$. Cette suite de Borels nous permet de
construire une d\'esingularisation de Demazure correspondant \`a
$\gb_1$ et $\gb'_1$ en prenant, si $\gb_x$ et $\gb_y$ sont deux termes
cons\'ecutifs de la suite le parabolique $\gb_x+\gb_y$ qui est minimal
mais diff\'erent d'un Borel. On construit ainsi pour tout $1\leq k\leq
n+1$ et tout $1\leq i\leq r_k-1$ les paraboliques
$\gp_{k,i}=\gb_{k,i}+\gb_{k,i+1}$, pour tout $1\leq k\leq n$ les
paraboliques $\gp_{k,r_k}=\gb_{k,r_k}+\gb_{k+1,1}$ et le parabolique
$\gp_{n+1,r_{n+1}}=\gb_{n+1,r_{n+1}}+\gb'_{n+1,r'_{n+1}}$ et de m\^eme
pour tout $1\leq k\leq n+1$ et tout $1\leq i\leq r'_k-1$ les paraboliques
$\gp'_{k,i}=\gb'_{k,i}+\gb'_{k,i+1}$ et pour tout $1\leq k\leq n$ les
paraboliques $\gp'_{k,r'_k}=\gb'_{k,r'_k}+\gb'_{k+1,1}$.
On pose alors
$Y'=\prod P'_{k,i}\times\prod P_{k,i}$ le
premier produit est effectu\'e dans l'ordre lexicographique et le
second dans l'ordre lexicographique invers\'e, on pose $G''=\prod
B'_{k,i}\times\prod B_{k,i}$ (avec les m\^emes ordres) et on pose
$X'=Y'/G''$. La d\'esingularisation de Demazure de
$\overline{B'_1B_1/B_1}$ est alors donn\'ee par $X'/B_1$ et un
morphisme $D$ de cette vari\'et\'e vers $G/B_1$. Le
morphisme $D$ est obtenu \`a partir de la multiplication de $Y'$
dans $G$. Or cette multiplication se factorise par $Y$ car les groupes
$P_{k,i}$ sont contenus dans $P_k$ et par passage au quotient on voit
que $D$ se factorise par $\pi$.

\vs 0.4 cm

\noi
{\bf Corollaire 1} : \textit{Le morphisme $\pi$ est une
 d\'esingularisation}

\vs 0.2 cm

\dm : 
On commence par montrer le cas des vari\'et\'es de Schubert :
$\overline{Bw'B/B}$ est exactement $\overline{w'P'_1B/B}$. On obtient
de cette fa\c con un ouvert (la $P'_1$-orbite $w'P'_1B/B$) lisse plus
grand que la cellule de Schubert. Le morphisme $\pi$ est un
isomorphisme au dessus de cet ouvert. En effet, cet ouvert contient la
cellule $Bw'B/B$ au dessus de laquelle la d\'esingularisation de
Demazure est un isomorphisme donc $\pi$ est aussi un isomorphisme au
dessus de cette cellule. De plus, le morphisme $\pi$ est invariant
sous l'action de $P_1$ donc l'orbite de la cellule lisse
pr\'ec\'edente (c'est exactement $w'P'_1B/B$) est une $P'_1$-orbite
lisse et $\pi$ est encore un isomorphisme au dessus de cette orbite.

Pour le cas g\'en\'eral on consid\`ere le diagramme commutatif suivant :
$$\begin{array}{ccc}
 X/B & \fl & \overline{Bw'B/B}\\
 \downarrow & & \downarrow\\
 X/P & \fl & \overline{w'P'_1P/P}
\end{array}$$
dont les fl\`eches verticales sont des fibrations en $P/B$ (pour la
seconde ceci vient du fait que comme $\gp_1+\gp'_1=\gb+\gb'=\gp+\gp'$,
la vari\'et\'e $\overline{Bw'B/B}$ est l'image r\'eciproque de
$\overline{w'P'_1B/B}$ par le morphisme de $G/B\fl G/P_1$)
ainsi comme le morphisme $\pi$ pour les Borels est une
d\'esingularisation et est donc birationnel, alors le morphisme $\pi$
de $X/P$ vers la vari\'et\'e de Schubert correspondante est
birationnel et c'est bien une d\'esingularisation.

\vs 0.4 cm

\noi
{\bf Remarque 2} :  La $P'_1$-orbite $w'P'_1B/B$ est le plus
grand ouvert lisse de $\overline{Bw'B/B}$ qui est une orbite
sous l'action d'un sous groupe de $G$ laissant stable
$\overline{Bw'B/B}$ (ceci vient du fait que $\gp'_1$ a \'et\'e
choisit maximal pour sa propri\'et\'e). Par ailleurs, pour tout groupe
$P$ contenant $B$ et contenu dans $P_1$, le morphisme
$\overline{Bw'P/P}\fl\overline{Bw'P_1/P_1}$ est une fibration en
$P_1/P$. Ainsi, les singularit\'es des vari\'et\'es
$\overline{Bw'P/P}$ pour de tels $P$ sont \textit{identiques} \`a
celle de $\overline{Bw'P_1/P_1}$. Les vari\'et\'es de Schubert du type
de $\overline{Bw'P_1/P_1}$ seront dits minimales : il n'existe pas de
parabolique $P$ contenant $P_1$ tel que la fl\`eche
$\overline{Bw'P_1/P_1}\fl\overline{Bw'P/P}$ est une fibration en
$P/P_1$.


\vs 0.4 cm

Notre construction nous permet de donner une condition n\'ecessaire et
suffisante pour qu'une vari\'et\'e de Schubert minimale soit
homog\`ene sous l'action d'un sous groupe de $G$ et on en d\'eduit une
condition suffisante de lissit\'e des vari\'et\'e de Schubert. Notons
$S$ la vari\'et\'e de Schubert $\overline{Bw'P/P}$ et $S_1$ la
vari\'et\'e de Schubert minimale associ\'ee $\overline{Bw'P_1/P_1}$ on
a alors le :

\vs 0.4 cm

\noi
{\bf Corollaire 2} : \textit{La vari\'et\'e de Schubert minimale $S_1$
est homog\`ene sous l'action d'un sous groupe de $G$ si et seulement si
$\gp_1\cap\gp'_1$ contient un Borel. Dans ce cas $S$ est lisse.}

\vs 0.2 cm

\dm :
Si $\gp_1\cap\gp'_1$ contient un Borel on a $X=P'_1\times^{P'_1\cap
P_1}P_1$. Si on quotiente par $P_1$ qui contient $P$ on a alors
$P'_1/(P'_1\cap P_1)=X/P_1\fl S_1$ est la d\'esingularisation de $S_1$
qui \'etait donc homog\`ene. Mais alors le morphisme $G/P\fl G/P_1$
nous donne le diagramme commutatif suivant :
$$\begin{array}{ccc}
 X/P & \fl & S\\
 \downarrow & & \downarrow\\
 X/P_1 & \fl & S_1
\end{array}$$
dont les fl\`eches verticales sont des fibrations en $P_1/P$ et donc,
comme la fl\`eche du bas est un isomorphisme, celle du haut est
\'egalement un isomorphisme. Ainsi on voit que $S$ est lisse car $X/P$
l'est.

R\'eciproquement, si $S_1$ est homog\`ene sous l'action d'un sous
groupe de $G$ alors le plus grand de ces sous groupes est $P'_1$ par
construction. Mais alors $P'_1/(P'_1\cap P_1)$ doit \^etre $S_1$ toute
enti\`ere. Ceci signifie que la cellule $P'_1/(P'_1\cap P_1)$ est
propre. Ceci n'est possible que si $P'_1\cap P_1$ contient un
Borel. Remarquons que pour que $S$ soit homog\`ene sous $P'_1$, il
faut et il suffit que $\P\cap P'_1$ contienne un Borel.

\vs 0.4 cm

\noi
{\bf Exemple 2} : On \'etudie notre d\'esingularisation pour certaines
vari\'et\'es de Schubert de $SL_4$. Soit $P_0\in L_0\subset H_0$ un
drapeau de $\p^3$. Consid\'erons les vari\'et\'es de Schubert
suivantes : $S=\{(P,L,H)\in\p^3\times\g\times{\check \p^3} / P\in
L\subset H \ {\rm{et}}\  {\rm{dim}}(L\cap L_0)\geq 1\}$,
$S'=\{\{(P,L,H)\in\p^3\times\g\times{\check \p^3} / P\in L\subset H,\
P\in L_0\ {\rm{et}}\ L\subset H_0\}$ et
$S''=\{\{(P,L,H)\in\p^3\times\g\times{\check \p^3} / P\in L\subset H,\
P\in H_0\ {\rm{et}}\ P_0\in H\}$, cette derni\`ere est la seule
vari\'et\'e de Schubert de $SL_4$ qui ne correspond pas \`a une
permutation vexillaire. Notre d\'esingularisation est alors donn\'ee
dans chacun des cas par
$X=\{(P,P',L,H,H')\in\p^3\times\p^3\times\g\times{\check
\p^3}\times{\check \p^3} / P\in L\subset H,\ P'\in L\subset H'\
{\rm{et}}\ \P'\in L_0\subset H'\}$, $X'=S'$ et
$X''=\{\{(P,L,L',H)\in\p^3\times\g\times\g\times{\check \p^3} / P\in
L\subset H,\ P\in L'\subset H\ {\rm{et}}\ P_0\in L'\subset H_0\}$. Ces
d\'esingularisations sont toutes bijectives sur le lieu lisse. Notons
que l'exemple de $S'$ montre que la condition suffisante de lissit\'e
du corollaire $2$ n'est pas n\'ecessaire.


\vs 0.4 cm

\noi
{\bf Remarque 3} : Si le lieu lisse d'une vari\'et\'e de Schubert $S$
est homog\`ene sous l'action d'un sous groupe de $G$ alors ce sous
groupe est $P_1$ et notre d\'esingularisation est bijective sur le
lieu lisse. Par ailleurs, on peut v\'erifier que si le groupe $G$ est
$SL_n$ pour $n\leq 4$, $SO_n$ pour $n\leq 6$ ou $Sp_n$ pour $n\leq 4$
alors notre d\'esingularisation est bijective sur le lieu lisse.

On peut encore affiner notre d\'esingularisation de la fa\c con
suivante : on remplace les paraboliques $P_i$ par un parabolique
contenant $P_i$ et contenu dans $w^{-1}\overline{BwB}$. Ceci ne change
pas les paraboliques $P_1$ et $P'_1$. Cette technique permet ainsi de
construire les d\'esingularisations
$X_1=\{(P,P',L,H)\in\p^3\times\p^3\times\g\times{\check \p^3} / P\in
L\subset H,\ P'\in L\ {\rm{et}}\ P'\in L_0\}$ et
$X_2=\{(P,L,H,H')\in\p^3\times\g\times{\check \p^3}\times{\check \p^3}
/ P\in L\subset H,\ L\subset H'\ {\rm{et}}\ L_0\subset H'\}$ de la
vari\'et\'e de Schubert $S$ de l'exemple $2$. On voit que cette
construction n'est plus canonique, il faut choisir un parabolique
contenant $P_i$ et contenu dans $w^{-1}\overline{BwB}$.

\vs 0.4 cm

\noi
{\bf{Remerciements}} : Je tiens ici \`a remercier mon directeur
de th\`ese Laurent Gruson pour toute l'aide qu'il m'a apport\'ee
durant la pr\'eparation de ce travail et \'egalement Patrick Polo pour
m'avoir signaler une erreur dans une premi\`ere version de cet article.

\section*{R\'ef\'erences}

\noi
[ACGH] Arbarello E. Cornalba M. Griffiths P.A. Harris J. : Geometry of
algebraic curves I, Grundlehren der Mathematischen Wissenschaften 267,
Springer Verlag, NewYork Berlin (1985).

\vs 0.1 cm

\noi
[B] Ballico E. : On the Hilbert scheme of curves in a smooth quadric,
Deformations of mathematical structures L\`od\`z/Lublin (1985/87)
Kluwer Acad. Publ. Dordrecht, 1989.

\vs 0.1 cm

\noindent
[BGG] Bernstein I.N. Gel'fand I.M. Gel'fand S.I. : Schubert cells, and
cohomology of the spaces $G/P$, Uspehi Mat. Nauk.(3) 68 (1973).

\vs 0.1 cm

\noi
[Bo] Borel A. : Linear algebraic groups, Second edition, Graduate Texts in Mathematics, 126. Springer-Verlag, New York (1991).

\vs 0.1 cm

\noindent
[D] Demazure M. : D\'esingularisation des vari\'et\'es de Schubert
g\'en\'eralis\'ees, Ann. Sci. ENS (4) 7 (1974).

\vs 0.1 cm

\noi
[DG] Demazure M. Gabriel P : Groupes alg\'ebriques, Tome I :
G\'eom\'etrie alg\'ebrique, g\'en\'eralit\'es, groupes
commutatifs. Masson \& Cie, Editeur, Paris; North-Holland Publishing
Co., Amsterdam (1970).

\vs 0.1 cm

\noindent
[FH] Fulton W. Harris J. : Representation theory, GTM 129 Springer
Verlag, New-York (1991).

\vs 0.1 cm

\noi
[Har] Harris J. : The genus of space curves, Math. Ann. 249 (1980) no. 3. 

\vs 0.1 cm

\noindent
[K1] Kempf G.R. : Vanishing theorems for flag maniflods,
Amer. j. math. 98 (1976).

\vs 0.1 cm

\noindent
[K2] Kempf G.R. : Linear systems, Ann. of Math. 103 (1976).

\vs 0.1 cm

\noindent
[Kl] Kleiman S.L. : The transversality of a general translate,
Compositio Math. 28 (1974).

\vs 0.1 cm

\noindent
[Ko] Koll\' ar J. : Rational curves on algebraic varieties, Ergebnisse
der Mathematik und ihrer Grenzgebiete 32, Springer Verlag, Berlin
(1996).

\vs 0.1 cm

\noindent
[KP] Kim B., Pandharipande R. : The connectedness of the moduli space
of maps to homogeneous spaces, preprint AG 0003168.

\vs 0.1 cm

\noindent
[LMS] Lakshmibai V., Musili C., Seshadri C.S. : Cohomology of line
bundles on $G/B$, Ann. Sci. ENS (4) 7 (1974).

\vs 0.1 cm

\noindent
[T] Thomsen J.F. : Irreducibility of $\overline{M}\sb
{0,n}(G/P,\beta)$, Internat. J. Math. 9 (1998) No. 3.

\end{document}